\documentclass[11pt,reqno]{amsart}

\usepackage[T1]{fontenc}
\usepackage[a4paper,margin=33mm]{geometry}

\usepackage{mathtools,bm}
\usepackage{amssymb}
\usepackage{microtype}
\usepackage{comment}
\usepackage[unicode,hidelinks]{hyperref}

\usepackage{float}
\usepackage{caption}
\usepackage{tikz}
\usetikzlibrary{calc,decorations.pathreplacing}
\usetikzlibrary{3d}

\usepackage{xpatch}
\makeatletter
\patchcmd{\@tocline}{\hfil}
{\nobreak\leaders\hbox{\ifnum#1<2\hfill\else$\m@th%
\mkern 4.5 mu\hbox{.}\mkern 4.5 mu$\fi}\hfill\nobreak}{}{}
\def\l@section{\@tocline{1}{10pt}{1pc}{}{\bfseries}}
\def\l@subsection{\@tocline{2}{0pt}{\dimexpr 1pc+2em}{}{}}
\makeatother

\theoremstyle{plain}
\newtheorem{theorem}{Theorem}[section]
\newtheorem{lemma}[theorem]{Lemma}
\newtheorem{corollary}[theorem]{Corollary}
\newtheorem{proposition}[theorem]{Proposition}

\theoremstyle{remark}
\newtheorem*{remark*}{Remark}

\numberwithin{equation}{section}

\newcommand{\eps}{\varepsilon}
\newcommand{\dx}{\mathrm{d} x}
\newcommand{\dy}{\mathrm{d} y}
\newcommand{\dz}{\mathrm{d} z}
\newcommand{\dw}{\mathrm{d} w}
\newcommand{\dt}{\mathrm{d} t}
\newcommand{\ds}{\mathrm{d} s}
\newcommand{\du}{\mathrm{d} u}
\newcommand{\dv}{\mathrm{d} v}
\newcommand{\rd}{\mathrm{d}}
\newcommand{\N}{\mathbb{N}}
\newcommand{\R}{\mathbb{R}}
\newcommand{\IR}{\int_{\R^d}}
\newcommand{\Sa}{S_\alpha}
\newcommand{\Da}{\mathcal D_\alpha}
\newcommand{\Ea}{\mathcal E_\alpha}
\newcommand{\Aa}{\mathcal A_\alpha}
\DeclareMathOperator{\tr}{tr}
\DeclareMathOperator{\supp}{supp}

\title[]{Lower bounds for Riesz interaction functionals}

\author[]{Nuno J. Alves}
\author[]{Athanasios E. Tzavaras}

\address[N. J. Alves and A. E. Tzavaras]{
      CEMSE Division, King Abdullah University of Science and Technology (KAUST), Thuwal 23955-6900, Saudi Arabia.}
\email{nuno.januarioalves@kaust.edu.sa}
\email{athanasios.tzavaras@kaust.edu.sa}

\begin{document}

\begin{abstract}
Motivated by the pressure formation problem for hydrodynamic models with nonlocal interactions, we classify strong Lebesgue lower bounds for the Riesz energy functional and for a determinant functional arising from the theory of compensated integrability. The classification covers all integrability exponents in $(0,\infty]$ within the class of confined log-concave densities. These geometric assumptions are motivated by counterexamples for general densities. On the whole space, neither functional controls any fixed Lebesgue quasi-norm. Under confinement, lower bounds fail for exponents greater than one, whereas $L^1$ bounds hold for the energy in every dimension and for the determinant only in dimension one. The log-concave estimates hold at sharp integrability exponents and follow from endpoint Lorentz bounds for confined quasiconcave densities and a weak-to-strong lemma. In particular, the energy estimates provide a geometric lower-bound counterpart to the classical Hardy--Littlewood--Sobolev inequality.
\end{abstract}

\maketitle
\thispagestyle{empty} 
\tableofcontents

\section{Introduction}

We consider the pressureless Euler--Riesz system with repulsive interactions, which models the motion of a compressible fluid driven by nonlocal forces. The equations take the form
\begin{equation}\label{eq:Euler--Riesz-system}
\begin{aligned}
\partial_t\rho &+\nabla\cdot(\rho u)=0,
\\[2mm]
\partial_t(\rho u) &+\nabla\cdot(\rho u\otimes u)
+\rho\nabla K_\alpha\ast\rho=0,
\end{aligned}
\end{equation}
and are posed in $(0,T) \times \R^d$, where $0<T<\infty$ and $d \in \N$ is the spatial dimension. Here, $\rho:[0,T)\times\R^d\to [0,\infty)$ denotes the fluid density and $u:[0,T)\times\R^d\to \R^d$ the velocity field. The kernel $K_\alpha$ describes the interaction between particles and is given by
\begin{equation}\label{eq:kernel_Ka}
K_\alpha(x)=\tfrac1{d-\alpha}|x|^{\alpha-d},
\qquad 0<\alpha<d.
\end{equation}
The associated Riesz interaction energy is
\begin{equation}\label{eq:E_functional}
\begin{split}
\Ea(\rho)
& =
\IR\rho\,(K_\alpha\ast\rho)\,\dx \\
& =
\tfrac1{d-\alpha}
\IR\IR
\rho(x)\,|x-y|^{\alpha-d}\,\rho(y)\,\dy\,\dx.
\end{split}
\end{equation}
For $d = 3$ and $\alpha=2$, $K_2$ is the Coulomb kernel, up to a constant factor, reducing the system to the repulsive pressureless Euler--Poisson system. Local well-posedness for the repulsive pressureless Euler--Riesz system is established in~\cite{choi2022wellposedness}.

For smooth solutions, the system~\eqref{eq:Euler--Riesz-system} formally conserves the total mass and energy,
\begin{equation}\label{eq:mass_energy_conservation}
\frac{\rd}{\dt}\IR\rho\,\dx=0,
\qquad
\frac{\rd}{\dt}
\bigg(
\IR\tfrac12\rho|u|^2\,\dx+\tfrac12\Ea(\rho)
\bigg)\!=0 , 
\end{equation}
which provide the basic a priori estimates for this problem. 

The goal of this work is to investigate whether one can extract further information from the conservation
laws~\eqref{eq:mass_energy_conservation}. To motivate this question, consider the Euler--Riesz system with a power-law pressure
\begin{equation}\label{eq:ERp}
\begin{aligned}
\partial_t\rho &+\nabla\cdot(\rho u) =0,
\\[2mm]
\partial_t(\rho u) &+\nabla\cdot(\rho u\otimes u)
+ \nabla \rho^\gamma + \rho\nabla K_\alpha\ast\rho =0,
\end{aligned}
\end{equation}
where $\gamma>1$. Since the kernel $K_\alpha$ is nonnegative, the energy identity
\begin{equation}
\frac{\rd}{\dt}
\bigg(
\IR\tfrac12\rho|u|^2+ \tfrac{1}{\gamma -1} \rho^\gamma + \tfrac12\rho\,(K_\alpha\ast\rho) \, \dx
\bigg)\!=0
\end{equation}
provides an a priori integrability estimate for the density, namely 
\[\rho \in L^{\infty}\big( (0,T); L^{\gamma}(\R^d)\big). \]
This bound can be improved using the theory of compensated
integrability (see~\cite[Section~3.2]{alves2026bilinear}),
which yields 
\[
\rho \in L^{\gamma+\frac1d}\big( (0,T); L^{\gamma+\frac1d}(\R^d)\big).
\]
The corresponding norms are controlled by the mass and energy at the initial time $t=0$, with bounds independent of $T$. Higher integrability of this kind is very useful in the existence theory for system~\eqref{eq:ERp}; see, for example,
\cite{carrillo2026global}, where such estimates and compensated compactness are employed to obtain existence of weak solutions under spherical symmetry.

Returning to the pressureless system~\eqref{eq:Euler--Riesz-system}, our objective is to investigate whether the interaction energy $\Ea(\rho)$ can control an $L^p$ norm of the density. Informally, the question is whether the energy provides a quantitative criterion for the formation of pressure. To place this problem in context, it is well known that kinetic descriptions of gases or dilute fluids lead to models involving nonlocal interactions. Examples are the Vlasov--Poisson system \cite{vlasov1968} for clouds of electrically charged particles, or various models for dilute fluids introduced by Irving and Kirkwood \cite{irving1950statistical}. In these models a quantitative description of the pressure is lacking even if intuitively the formation of pressure is expected.

We study this question in the hydrodynamic setting of the pressureless Euler--Riesz system. To this end, we follow two approaches, based on lower bounds for (i) the energy $\Ea(\rho)$ defined in \eqref{eq:E_functional} and (ii) the functional $\Da(\rho)$ defined in~\eqref{eq:D_functional}, which is motivated by compensated integrability and whose time integral is controlled by the initial mass and total energy, as outlined below.

To throw light on (ii), let us review  Serre's theory of compensated integrability~\cite{serre2018divergence,serre2019compensated},
which combines the divergence-free condition with symmetry and positivity of a stress tensor to obtain determinant estimates. For further developments on the subject, the reader is referred to \cite{derosa2020upper, derosa2023fine, guerra2024compensation}.
To apply compensated integrability to~\eqref{eq:Euler--Riesz-system}, we introduce the Riesz interaction stress 
tensor \begin{equation}\label{eq:Sa}
\Sa(\rho)(x)
=
\tfrac12\int_0^1\!\IR
\rho\big(x+(\theta-1)y\big)
\rho\big(x+\theta y\big)
|y|^{\alpha-d-2}\,y\otimes y\,\dy\,\rd\theta.
\end{equation}
This tensor is symmetric, positive semidefinite and satisfies (see \cite{alves2026bilinear})
\begin{equation}\label{eq:Ka_Sa_identity}
\rho\,\nabla K_\alpha\ast\rho
=
\nabla\cdot\Sa(\rho).
\end{equation}
Integrating its trace gives the identity
\begin{equation}\label{eq:trace-stress-energy}
\Ea(\rho)
=
\tfrac2{d-\alpha}\IR\operatorname{tr}\Sa(\rho)\,\dx,
\end{equation}
so the Riesz energy is (up to a constant) the integral of the trace of the interaction stress.
Using~\eqref{eq:Ka_Sa_identity}, the system~\eqref{eq:Euler--Riesz-system} may be written in a (space-time) divergence-free form,
\begin{equation}\label{eq:Euler--Riesz-divergence}
\nabla_{t,x}\cdot\Aa(\rho,u)=0,
\end{equation}
where the matrix 
\begin{equation}\label{eq:Euler--Riesz-tensor}
\Aa(\rho,u)
:=
\begin{pmatrix}
\rho&(\rho u)^\top\\
\rho u&\rho u\otimes u+\Sa(\rho)
\end{pmatrix}
\end{equation}
is symmetric and positive semidefinite. A direct computation gives
\begin{equation}\label{eq:det-A}
\det\Aa(\rho,u)=\rho\det\Sa(\rho).
\end{equation}
The theory of compensated integrability then provides an estimate 
\[
\int_0^T\!\IR
\big(\!\det\Aa(\rho,u)\big)^{1/d}\,\dx\,\dt
=
\int_0^T\Da(\rho)\,\dt  \le C
\]
in terms of the initial mass and total energy; see \cite{serre2018divergence, serre2019compensated},
and \cite{alves2026bilinear} in the context of \eqref{eq:Euler--Riesz-system}.
This suggests seeking lower bounds for the functional
\begin{equation}\label{eq:D_functional}
\Da(\rho)
:=
\IR
\big(\rho\,\det\Sa(\rho)\big)^{1/d}\,\dx.
\end{equation}

There are alternative ways to represent the stress. For instance, the approach in~\cite{GLT2017} for Euler--Poisson systems can be adapted to Euler--Riesz systems using the Caffarelli--Silvestre extension. This representation is used in~\cite{chaudhuri2026generalized} to formulate a notion of dissipative solution and prove weak--strong uniqueness for the repulsive barotropic Euler--Riesz system. The associated stress is symmetric but, unlike $\Sa(\rho)$, is not positive semidefinite.

The interaction stress $\Sa(\rho)$ has the structure of bilinear fractional integrals. Uniform estimates for the associated operators are obtained in~\cite{alves2026bilinear, alves2026uniform}. The present work concerns the complementary problem of obtaining lower bounds for the functionals $\Ea$ and $\Da$ in terms of Lebesgue quasi-norms of the density. This leads to a constrained minimization problem:  given a class $\mathcal C$ of nonnegative densities, one considers the infima
\[
\inf_{\substack{\rho\in\mathcal C\\
\|\rho\|_{L^p(\R^d)}=1}}
\Ea(\rho),
\qquad
\inf_{\substack{\rho\in\mathcal C\\
\|\rho\|_{L^q(\R^d)}=1}}
\Da(\rho).
\]
The objective then becomes to identify exponents $p,q\in(0,\infty]$ and the geometric assumptions on~$\mathcal C$ under which these infima are strictly positive. We begin without geometric restrictions and then examine the effect of confinement, quasiconcavity, and log-concavity.
For the classes considered here, positivity of the corresponding infima is equivalent to uniform lower bounds of the form
\begin{equation}\label{eq:introlbineq}
\Ea(\rho)\geq C\|\rho\|_{L^p(\R^d)}^2
\qquad\text{and}\qquad
\Da(\rho)\geq C\|\rho\|_{L^q(\R^d)}^\beta,
\end{equation}
where
\begin{equation}\label{eq:beta}
\beta:=2+\frac1d.
\end{equation}
The inequalities \eqref{eq:introlbineq} are preserved by the transformation $\rho \mapsto \lambda \rho$, $\lambda > 0$, as the functionals $\Ea$ and $\Da$ are homogeneous of degrees $2$ and $\beta$, respectively.

The selection of norms in  these lower bounds turns out to depend on geometric assumptions  on the density. Without such assumptions, neither functional controls any fixed Lebesgue quasi-norm on the whole space; see Theorem~\ref{thm:whole-space-negative}. For arbitrary densities supported in a fixed ball, the energy admits an $L^1$ lower bound in every dimension, while the corresponding determinant bound holds only in dimension one; see Propositions~\ref{prop:E_L1_confined}
and~\ref{prop:D-confined-L1-endpoint}
and Theorem~\ref{thm:confined-negative}.

Stronger lower bounds are achieved when confinement is combined with log-concavity of the density. Under these assumptions, strong Lebesgue estimates are obtained at explicit critical exponents; see Theorem~\ref{thm:ED-log-concave}. Their proof uses endpoint Lorentz bounds for quasiconcave densities and a weak-to-strong lemma. Together with the sharpness theorem and finite-measure embeddings, these estimates provide a complete classification of strong Lebesgue lower bounds over all exponents in $(0,\infty]$ within the class of confined log-concave densities. For the energy functional, they comprise a geometric lower-bound counterpart to the classical Hardy--Littlewood--Sobolev inequality.

Turning to the evolution problem, estimates with Lebesgue exponents greater than one are of particular interest, since they would provide integrability beyond conservation of mass. Applying such lower bounds along the flow would require
that the relevant geometric assumptions be preserved. If we suppose that on a time interval $(0,T_*)$ the density belongs, for almost all times, to a class for which the relevant lower bound holds, then energy conservation and compensated integrability would yield, respectively,
\[
\rho\in L^\infty\bigl((0,T_*);L^p(\R^d)\bigr),
\qquad
\rho\in L^\beta\bigl((0,T_*);L^q(\R^d)\bigr).
\]
The present work establishes the underlying static inequalities; it does not address the challenging question of determining under which additional conditions the geometric assumptions might be preserved by the evolution. The functionals $\Ea(\rho)$ and $\Da(\rho)$ describe, respectively, the trace and determinant mechanisms through which the interaction stress may supply pressure-like integrability for densities with a specified geometric structure. The main results are stated precisely in Section \ref{sec:main-results}.

\section{Statement of main results} \label{sec:main-results}

We begin by fixing the notation and terminology used throughout the paper. By a \emph{density} we mean a nonnegative measurable function on $\R^d$. For $m\in\N$, $x\in\R^m$, and $r>0$, we denote by $B_r^m(x)$ the open ball in $\R^m$ centered at $x$ with radius $r$, and write $B_r^m:=B_r^m(0)$. When $m=d$, we omit the superscript and simply write $B_r(x)$ and $B_r$. We say that a density $\rho$ is supported in $B_R$ if $\rho=0$ almost everywhere in $\R^d\setminus B_R$. For a measurable set
$E\subseteq\R^d$, we denote its Lebesgue measure by $|E|$.

For $0<p\leq\infty$, we use the usual Lebesgue norms and
quasi-norms
\[
\|\rho\|_{L^p(\R^d)}
:=
\begin{dcases}
\displaystyle
\bigg(\IR \rho(x)^p\,\dx\bigg)^{1/p},
&\qquad 0<p<\infty,\\[4mm]
\operatorname*{ess\,sup}_{x\in\R^d}\rho(x),
&\qquad p=\infty.
\end{dcases}
\]
When $\rho$ is supported in $B_R$, these quantities agree with
$\|\rho\|_{L^p(B_R)}$.

For real symmetric $d\times d$ matrices $X$ and $Y$, we write $X\succeq Y$ if $X-Y$ is positive semidefinite; this is the L\"{o}wner order. We shall use that the determinant is monotone on the positive semidefinite cone with respect to this order: if $X\succeq Y\succeq0$, then $\det X\geq\det Y$; see~\cite[Corollary~7.7.4]{horn2013matrix}. The $d\times d$ identity matrix is denoted by $I_d$.
The symbols $c$ and $C$, possibly with subscripts indicating the parameters on which they may depend, denote positive constants whose values may change from line to line.

With these conventions in place, we state the main results. Section~\ref{subsec:general-densities} concerns lower bounds without any assumptions on the shape of the density, both on the whole space and for densities supported in a fixed ball. Then, in Section~\ref{sec:logconcave} we present lower bounds for confined log-concave densities and a theorem regarding their sharpness.

\subsection{Lower bounds for general densities}
\label{subsec:general-densities}

We first consider the whole space and construct a counterexample showing that lower bounds cannot hold for general densities. The example is based on dispersion: a set of fixed measure is divided into an increasing number of smaller components placed
farther apart. This yields a single sequence $\{\rho_n\}_{n \in \N}$ of densities whose $L^q$ quasi-norms are independent of $n$ for every $0<q\leq\infty$, while both $\mathcal E_\alpha(\rho_n)$ and $\Da(\rho_n)$ tend to zero.

\begin{theorem}\label{thm:whole-space-negative}
Let $d\in\N$ and $0<\alpha<d$. There exists a sequence $\{\rho_n\}_{n\in\N}$ of densities such that
\[
\|\rho_n\|_{L^q(\R^d)}=c_{d,q}>0
\qquad\text{for all }0< q\le\infty,
\]
with constants $c_{d,q}$ independent of $n$, while
\begin{equation} \label{eq:whole-space-negative}
\mathcal E_\alpha(\rho_n)\to0, \qquad \Da(\rho_n)\to 0,
\qquad\text{as }n\to\infty.
\end{equation}
\end{theorem}

Next, we consider densities supported in a fixed bounded domain. We again produce a counterexample, now based on a generalization of the so-called Dirac comb. Confinement prevents arbitrarily large separations but still permits fragmentation into many small components. For both $\mathcal E_\alpha$ and $\mathcal D_\alpha$, this confined
fragmentation rules out a uniform lower bound in terms of the $L^q$-norm for every $1<q<\infty$. For the determinant functional, it also rules out the $L^1$ endpoint in dimensions $d\geq2$.

\begin{theorem}\label{thm:confined-negative}
Let $d\in\N$, $0<\alpha<d$, and $1\leq q<\infty$.
There exists a sequence $\{\rho_n\}_{n\in\N}$ of densities supported in $[0,1]^d$ such that
\[
\|\rho_n\|_{L^q([0,1]^d)}=c_{d,q}>0,
\]
with $c_{d,q}$ independent of $n$, and
\begin{enumerate}
\item[(i)]
if $q>1$, then
\begin{equation}\label{eq:E-confined-negative}
\Ea(\rho_n)\to0
\qquad\text{as }n\to\infty,
\end{equation}

\item[(ii)]
if $d=1$ and $q>1$, or if $d\geq2$ and $q\geq1$, then
\begin{equation}\label{eq:D-confined-negative}
\Da(\rho_n)\to0
\qquad\text{as }n\to\infty.
\end{equation}
\end{enumerate}
\end{theorem}

The case $q=\infty$ is immediate. Indeed, for $0<\eps<R$, set $\rho_\eps:=\chi_{B_\eps}$. Then $\rho_\eps$ is supported in $B_R$,
$\|\rho_\eps\|_{L^\infty(\R^d)}=1$, and
\[
\Ea(\rho_\eps)
=
\eps^{d+\alpha}\Ea(\chi_{B_1}),
\qquad
\Da(\rho_\eps)
=
\eps^{d+\alpha}\Da(\chi_{B_1}).
\]
Hence both functionals tend to zero as $\eps\to0$, and no uniform lower bound can hold when $q=\infty$.

\begin{remark*} 
In the whole-space and in the confined space constructions, the densities are constant multiples of characteristic functions of suitable unions of disjoint sets. The same conclusions hold with densities in
$C_c^\infty(\R^d)$. Indeed, on each component one may replace the characteristic function by a translated and rescaled copy of a fixed
nonnegative smooth bump, while keeping the same amplitude. The support conditions are preserved and the relevant $L^q$-norms remain independent of $n$.
\end{remark*}

At the $L^1$ endpoint, confinement provides a positive result. Every density $\rho$ supported in $B_R$ satisfies
\[
\Ea(\rho)
\geq
C_{\alpha,d,R}\,
\|\rho\|_{L^1(B_R)}^2,
\qquad 0<\alpha<d,
\]
in every dimension; see Proposition~\ref{prop:E_L1_confined}.
For the determinant functional, we prove the corresponding estimate in dimension one:
\[
\Da(\rho)
\geq
C_{\alpha,R}\,
\|\rho\|_{L^1(B_R)}^3,
\qquad d=1,\; \;  0<\alpha<1,
\]
for  $\rho \ge 0$ supported in $B_R\subseteq\R$;
see Proposition~\ref{prop:D-confined-L1-endpoint}.
Neither estimate requires geometric assumptions beyond confinement.

\subsection{Sharp lower bounds for log-concave densities}\label{sec:logconcave}

The counterexamples in Theorems~\ref{thm:whole-space-negative}
and~\ref{thm:confined-negative} are based on fragmentation. To obtain lower bounds beyond those supplied by confinement alone, we impose geometric conditions on the class of densities considered: we work with either quasiconcave or log-concave densities.

Recall that a density $\rho$ is called \emph{log-concave} if
\begin{equation}\label{eq:log-concavity}
\rho\big((1-\theta)x+\theta y\big)
\geq
\rho(x)^{1-\theta}\rho(y)^\theta
\end{equation}
for every $x,y\in\R^d$ and every $0<\theta<1$.
Also, a density is called \emph{quasiconcave} if all its superlevel sets are convex. These notions are well-developed in convex geometry and analysis; see, for instance, \cite{brascamp1976extensions, bobkov1999isoperimetric,
schneider2014convex, colesanti2017minkowski, haddad2025affine}.
The log-concave class is a subclass of the quasiconcave class, as every positive superlevel set of a log-concave density is convex. 

For $0<\alpha<d$, define the critical exponents
\begin{equation}\label{eq:p_alpha}
p_\alpha
:=
\begin{dcases}
\displaystyle
\frac{2}{\alpha+1},
&\qquad 0<\alpha\leq1,\\[2mm]
1,
&\qquad 1\leq\alpha<d,
\end{dcases}
\end{equation}
and
\begin{equation}\label{eq:q_alpha}
q_\alpha
:=
\begin{dcases}
\displaystyle
\frac{2+1/d}{\alpha+1},
&\qquad 0<\alpha\leq1,\\[3mm]
\displaystyle
\frac{2d+1}{2d+\alpha-1},
&\qquad 1\leq\alpha\leq3,\\[3mm]
\displaystyle
\frac{2d+1}{2d+2},
&\qquad 3\leq\alpha<d.
\end{dcases}
\end{equation}
The expressions agree at the common endpoints.
In particular,
\[
q_\alpha\geq1
\qquad\text{if and only if}\qquad
\alpha\leq2.
\]

The following theorem, proved in
Section~\ref{sec:lorentz-quasiconcave}, gives strong Lebesgue lower bounds at these exponents.

\begin{theorem}
\label{thm:ED-log-concave}
Let $d\in\N$, $0<\alpha<d$, and $R>0$. There exists a positive constant
$C=C(\alpha,d,R)$ such that
\begin{equation}\label{eq:ED-log-concave}
\mathcal E_\alpha(\rho)
\ge
C\,
\|\rho\|_{L^{p_\alpha}(B_R)}^2
\qquad
\text{and}
\qquad
\Da(\rho)
\ge
C\,
\|\rho\|_{L^{q_\alpha}(B_R)}^\beta
\end{equation}
for all log-concave densities $\rho$ supported in $B_R$.
\end{theorem}

Since $|B_R|<\infty$, H\"older's inequality gives
\[
\|\rho\|_{L^a(B_R)}
\leq
|B_R|^{1/a-1/b}\|\rho\|_{L^b(B_R)},
\qquad 0<a<b<\infty.
\]
Thus the corresponding lower bounds hold for every
$0<p<p_\alpha$ and $0<q<q_\alpha$, respectively.

The following theorem shows that the critical exponents
$p_\alpha$ and $q_\alpha$ cannot be increased.

\begin{theorem}\label{thm:sharpness_1}
Let $d\in\N$, $0<\alpha<d$, and $0<r<\infty$. There exists a family
$\{\rho_{\eps,r}\}_{\eps\in(0,1)}$ of log-concave densities such that
\[
\supp\rho_{\eps,r}\subseteq B_{\sqrt2},
\qquad
\|\rho_{\eps,r}\|_{L^r(\R^d)}=1,
\]
and
\begin{enumerate}
\item[(i)]
if $r>p_\alpha$, then
\[
\Ea(\rho_{\eps,r})\to0
\qquad
\text{as }\eps\to0.
\]

\item[(ii)]
if $r>q_\alpha$, then
\[
\Da(\rho_{\eps,r})\to0
\qquad
\text{as }\eps\to0.
\]
\end{enumerate}
\end{theorem}

The construction uses normalized characteristic functions of thin plates.
A simple dilation argument places these counterexamples in any prescribed ball $B_R$. Together with Theorem~\ref{thm:ED-log-concave}, the finite-measure embedding above, and the shrinking-ball construction for $r=\infty$, this gives the exact ranges
\[
0<p\leq p_\alpha
\qquad\text{and}\qquad
0<q\leq q_\alpha
\]
for strong Lebesgue lower bounds for $\Ea$ and $\Da$,
respectively, uniformly over log-concave densities supported
in $B_R$.

\subsection{Relation to HLS inequalities}

The Riesz energy $\Ea$ is the quadratic form associated with the Riesz potential of order $\alpha$. The classical Hardy--Littlewood--Sobolev inequality gives
\[
\Ea(\rho)
\leq
C_{\alpha,d} \,
\|\rho\|_{L^{\frac{2d}{d+\alpha}}(\R^d)}^2
\]
for every nonnegative $\rho\in L^{\frac{2d}{d+\alpha}}(\R^d)$. Its sharp form, including the optimal constant and the characterization of the extremizers, was established by Lieb~\cite{lieb1983sharp}; see also~\cite{frank2010inversion}. More recently, Haddad and Ludwig~\cite{haddad2025affine} obtained refinements that strengthen the classical inequality, together with reverse inequalities for the
associated affine functional under log-concavity.

The estimates proved here concern inequalities in the opposite direction. For log-concave densities supported in $B_R$, Theorem~\ref{thm:ED-log-concave} provides a lower bound with the critical exponent $p_\alpha$. This is not, in general, a literal reversal of the classical Hardy--Littlewood--Sobolev inequality. The exponent $2d/(d+\alpha)$ in the classical upper bound is dictated by whole-space scaling, whereas the exponent $p_\alpha$ in the lower bound arises from geometric estimates for bounded convex sets. In particular,
\[
p_\alpha
<
\frac{2d}{d+\alpha}
\]
when $d\geq2$, while the two exponents agree when $d=1$.

In dimension one, the energy and the squared $L^{2/(1+\alpha)}$ norm have the same spatial scaling. Rescaling the lower bound of Theorem~\ref{thm:ED-log-concave} from the unit interval gives a constant independent of $R$. Since the truncations $\rho\,\chi_{(-R,R)}$ remain log-concave, monotone convergence as $R\to\infty$ removes the confinement assumption. Together with the classical Hardy--Littlewood--Sobolev inequality, this gives the following result.

\begin{proposition}\label{prop:one-dimensional-energy-equivalence}
Let $d=1$ and $0<\alpha<1$. There exist positive constants
$c_\alpha$ and $C_\alpha$ such that
\begin{equation}\label{eq:one-dimensional-energy-equivalence}
c_\alpha \,
\|\rho\|_{L^{\frac{2}{1+\alpha}}(\R)}^2
\leq
\Ea(\rho)
\leq
C_\alpha \,
\|\rho\|_{L^{\frac{2}{1+\alpha}}(\R)}^2
\end{equation}
for all log-concave densities $\rho\in L^{\frac{2}{1+\alpha}}(\R)$.
\end{proposition}

The one-dimensional lower bound can also be deduced from the reverse inequality in~\cite[Theorem~1.2]{haddad2025affine}, since in dimension one the associated affine functional reduces, up to a constant, to the Riesz energy.

Classical rearrangement theory provides a further comparison.
If $\rho^\ast$ denotes the symmetric decreasing rearrangement of $\rho$,
then the Riesz rearrangement inequality~\cite{brascamp1974rearrangement},
with equality cases studied in~\cite{burchard1996riesz}, gives
\[
\Ea(\rho)\leq\Ea(\rho^\ast).
\]
Since rearrangement preserves Lebesgue quasi-norms but can increase the energy, this approach does not directly reduce the present lower-bound problem to symmetric decreasing densities.

We finally distinguish the present results from the reverse
Hardy--Littlewood--Sobolev inequalities of~\cite{carrillo2019reversehls}. Those inequalities concern lower bounds for interaction energies with positive-power kernels $|x-y|^\lambda$, $\lambda>0$, whereas we consider the singular decreasing kernel $|x-y|^{\alpha-d}$, $0<\alpha<d$.

\subsection{Outline of the proofs}

In Section~\ref{sec:whole-space-obstructions}, we prove Theorem~\ref{thm:whole-space-negative} by placing many small disjoint components increasingly far apart. In Section~\ref{sec:confined-obstructions}, Theorem~\ref{thm:confined-negative} is proved using densities supported on many small disjoint components inside a fixed cube.
In Section~\ref{sec:lower-bounds-conf-alone}, we establish the $L^1$ lower bounds of
Propositions~\ref{prop:E_L1_confined}
and~\ref{prop:D-confined-L1-endpoint}
that hold under confinement alone.

We then restrict attention to densities satisfying geometric conditions. Section~\ref{sec:box-estimates} provides lower bounds for characteristic functions of boxes. In Section~\ref{sec:lorentz-quasiconcave}, we use these bounds
to prove Theorem~\ref{thm:ED-log-concave}.
John's ellipsoid theorem is employed to extend the box estimates to characteristic functions of
convex sets. Applying these estimates to the superlevel sets of quasiconcave densities and then integrating over the level parameter yields
the endpoint Lorentz bounds of Theorem~\ref{thm:ED-quasiconcave}. The corresponding weak Lebesgue estimates follow immediately. We then apply a weak-to-strong lemma for log-concave densities,
Lemma~\ref{lem:log-concave-weak-strong}, to obtain the strong Lebesgue estimates.
Finally, Section~\ref{sec:sharpness} proves Theorem~\ref{thm:sharpness_1} using thin plates to establish the optimality of the critical exponents.

\section{Densities spreading in space}\label{sec:whole-space-obstructions}

In this section we prove Theorem~\ref{thm:whole-space-negative}. The two conclusions follow from the same estimate for the auxiliary functional
\begin{equation}\label{eq:J_functional}
J_\alpha(f,g)
:=
\IR \IR
f(z) \, g(w) \, |z-w|^{\alpha-d}\,\dz\,\dw,
\qquad
J_\alpha(f):=J_\alpha(f,f),
\end{equation}
defined for nonnegative measurable functions $f$ and $g$ on $\R^d$.

The example consists of $n$ balls of radius $n^{-1/d}$ placed far apart. Their total measure is independent of $n$, while the next lemma shows that the corresponding value of $J_\alpha$ tends to zero.
For $n\in\N$, set
\begin{equation}\label{eq:whole-space-parameters}
r_n:=n^{-1/d},
\qquad
x_{n,k}:=4kn\,e_1,
\qquad
k=1,\dots,n,
\end{equation}
where $e_1=(1,0,\dots,0) \in \R^d$, define
\begin{equation}\label{eq:whole_space_set}
E_n
:=
\bigcup_{k=1}^n B_{r_n}(x_{n,k}) ,
\end{equation}
and consider the sequence of densities $\rho_n = \chi_{E_n}$.

\begin{figure}[H]
\centering

\begin{minipage}[b]{0.47\textwidth}
\centering
\begin{tikzpicture}[scale=0.82]

    \draw[->,thick,gray!70] (0,0) -- (6.8,0) node[below]{};
    \draw[->,thick,gray!70] (0,-1.2) -- (0,1.2) node[left]{};

    \fill[gray!75] (1.9,0) circle (0.36);
    \draw[gray!85] (1.9,0) circle (0.36);

    \fill[gray!75] (4.9,0) circle (0.36);
    \draw[gray!85] (4.9,0) circle (0.36);

\end{tikzpicture}

\par\smallskip
{\small (a) The set $E_n$ for $n=2$.}
\end{minipage}
\hspace{0.035\textwidth}
\begin{minipage}[b]{0.47\textwidth}
\centering
\begin{tikzpicture}[scale=0.82]

    \draw[->,thick,gray!70] (0,0) -- (6.8,0) node[below]{};
    \draw[->,thick,gray!70] (0,-1.2) -- (0,1.2) node[left]{};

    \fill[gray!75] (1.2,0) circle (0.22);
    \draw[gray!85] (1.2,0) circle (0.22);

    \fill[gray!75] (2.7,0) circle (0.22);
    \draw[gray!85] (2.7,0) circle (0.22);

    \fill[gray!75] (4.2,0) circle (0.22);
    \draw[gray!85] (4.2,0) circle (0.22);

    \fill[gray!75] (5.7,0) circle (0.22);
    \draw[gray!85] (5.7,0) circle (0.22);

\end{tikzpicture}

\par\smallskip
{\small (b) The set $E_n$ for $n=4$.}
\end{minipage}

\caption{Schematic illustration of $E_2$ and $E_4$ in dimension $d=2$.
As~$n$ increases, the~$n$ disjoint balls become smaller and their centers farther apart. The drawings are not to scale: the minimum distance between distinct centers is~$8$ for $E_2$ and $16$ for $E_4$.}
\label{fig:whole-space-construction}
\end{figure}
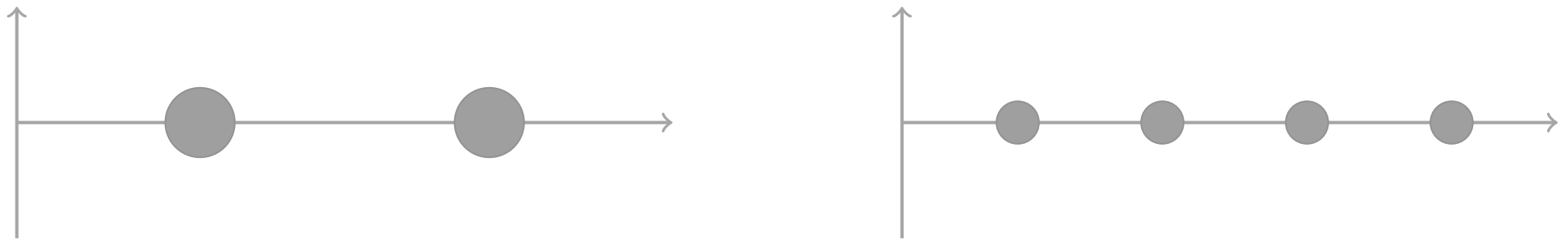

\begin{lemma}
\label{lem:whole-space-J-estimate}
Let $d\in\N$ and $0<\alpha<d$. For every $n\in\N$, the balls appearing in~\eqref{eq:whole_space_set} are pairwise disjoint,
\[
\|\chi_{E_n}\|_{L^q(\R^d)}
=
\begin{dcases}
|B_1|^{1/q},
\qquad &0< q<\infty,\\[1mm]
1,
\qquad &q=\infty,
\end{dcases}
\]
and
\begin{equation}\label{eq:whole-space-J-estimate}
J_\alpha(\chi_{E_n})
\le
c_{\alpha,d}
\big(
n^{-\alpha/d}
+
n^{\alpha-d}
\big).
\end{equation}
In particular, $J_\alpha(\chi_{E_n})\to0 \;  \text{as }n\to\infty$.
\end{lemma}

\begin{proof}
Since $r_n\le1$ and
\[
|x_{n,i}-x_{n,j}|
=
4n|i-j|
\ge
4n \geq 4
\qquad
(i\neq j),
\]
the balls $B_{r_n}(x_{n,k})$ are pairwise disjoint. Moreover,
\[
|E_n|
=
n \,|B_{r_n}|
=
n \,r_n^d |B_1|
=
|B_1|,
\] 
\[
\| \chi_{E_n} \|^q_{L^q(\R^d)} = |E_n| = |B_1|  \quad \text{for } q < \infty,
\]
while $\|\chi_{E_n}\|_{L^\infty(\R^d)}=1$.

Set now $B_{n,k}:=B_{r_n}(x_{n,k})$, $k=1,\dots,n$ and observe that, 
by bilinearity,
\[
J_\alpha(\chi_{E_n})
=
\sum_{k=1}^nJ_\alpha(\chi_{B_{n,k}})
+
\sum_{\substack{1\le i,j\le n\\i\neq j}}
J_\alpha(\chi_{B_{n,i}},\chi_{B_{n,j}}).
\]
We estimate the first sum. For each $k$, the change of variables
\[
z=x_{n,k}+r_nu,
\qquad
w=x_{n,k}+r_nv
\]
gives
\begin{align*}
J_\alpha(\chi_{B_{n,k}})
&=
r_n^{2d} \, r_n^{\alpha-d}
\int_{B_1} \int_{B_1}
|u-v|^{\alpha-d}\,\du\,\dv
\\
&=
r_n^{d+\alpha}J_\alpha(\chi_{B_1}).
\end{align*}
The quantity $J_\alpha(\chi_{B_1})$ is finite because the kernel
$(u,v) \mapsto |u-v|^{\alpha-d}$ is locally integrable in $\R^{2d}$. Since $r_n=n^{-1/d}$,
\begin{align*}
\sum_{k=1}^nJ_\alpha(\chi_{B_{n,k}})
&=
n \,r_n^{d+\alpha}J_\alpha(\chi_{B_1})
\\
&=
n^{-\alpha/d}J_\alpha(\chi_{B_1}).
\end{align*}

Next, we estimate the second sum. If $i\neq j$,
$z\in B_{n,i}$ and $w\in B_{n,j}$, then
\begin{align*}
|z-w|
&\ge
|x_{n,i}-x_{n,j}|-2r_n
\\
&\ge
4n-2
\ge
2n.
\end{align*}
Since $\alpha-d<0$, we have
\[
|z-w|^{\alpha-d}
\le
(2n)^{\alpha-d} , 
\]
\[
|B_{n,k}|
=
r_n^d|B_1|
=
\frac{|B_1|}{n},
\]
and 
\begin{align*}
\sum_{\substack{1\le i,j\le n\\i\neq j}}
J_\alpha(\chi_{B_{n,i}},\chi_{B_{n,j}})
&\le
(2n)^{\alpha-d}
\sum_{\substack{1\le i,j\le n\\i\neq j}}
|B_{n,i}|\,|B_{n,j}|
\\
&=
(2n)^{\alpha-d}
\frac{n(n-1)}{n^2}|B_1|^2
\\
&\le
c_{\alpha,d} \, n^{\alpha-d}.
\end{align*}
Combining the two estimates, we arrive at
\begin{align*}
J_\alpha(\chi_{E_n})
\le
n^{-\alpha/d}J_\alpha(\chi_{B_1})
+
c_{\alpha,d} \,n^{\alpha-d}
\le
c_{\alpha,d}
\big(
n^{-\alpha/d}
+
n^{\alpha-d}
\big),
\end{align*}
which proves~\eqref{eq:whole-space-J-estimate}.
\end{proof}

\begin{proof}[Proof of Theorem~\ref{thm:whole-space-negative}]
Let $\rho_n := \chi_{E_n}$ with $E_n$ defined by~\eqref{eq:whole_space_set}. The convergence
\[
J_\alpha(\rho_n)\to0
\qquad
\text{as }n\to\infty
\]
follows from Lemma~\ref{lem:whole-space-J-estimate}, and using that 
$\mathcal E_\alpha(\rho_n)
=
\tfrac1{d-\alpha}J_\alpha(\rho_n),
$
we conclude
\[
\mathcal E_\alpha(\rho_n)\to0
\qquad
\text{as }n\to\infty.
\]

It remains to estimate the determinant functional. Since $\Sa(\rho_n)(x)$ is positive
semidefinite, the arithmetic-geometric mean inequality gives
\[
\big(\!\det\Sa(\rho_n)(x)\big)^{1/d}
\le
\frac1d\tr\Sa(\rho_n)(x)
\]
and thus
\[
\Da(\rho_n)
\le
\frac1d
\int_{\R^d}
\rho_n(x)^{1/d}\tr\Sa(\rho_n)(x)\,\dx.
\]
Since $0\le\rho_n^{1/d}\le1$, 
\[
\Da(\rho_n)
\le
\frac1d
\int_{\R^d}
\tr\Sa(\rho_n)(x)\,\dx.
\]

Taking the trace of $\Sa(\rho_n)(x)$ and using Tonelli's theorem, we obtain
\begin{align*}
\int_{\R^d}\tr\Sa(\rho_n)(x)\,\dx
&=
\frac12
\int_0^1
\int_{\R^d}
\int_{\R^d}
\rho_n\big(x+(\theta-1)y\big)\,
\rho_n\big(x+\theta y\big)\,
|y|^{\alpha-d}
\,\dy\,\dx\,\rd\theta.
\end{align*}
For fixed $\theta$, the change of variables
\[
z=x+(\theta-1)y,
\qquad
w=x+\theta y,
\]
has Jacobian one and satisfies $y=w-z$, thus implying
\[
\int_{\R^d}\tr\Sa(\rho_n)(x)\,\dx
=
\frac12
\IR \IR
\rho_n(z) \, \rho_n(w) \,|z-w|^{\alpha-d}
\,\dz\,\dw
=
\frac12J_\alpha(\rho_n).
\]
Thus
\[
\Da(\rho_n)
\le
\frac1{2d}J_\alpha(\rho_n)
\]
and the proof follows from Lemma~\ref{lem:whole-space-J-estimate}.
\end{proof}

\section{Densities with confined support} \label{sec:confined-obstructions}

In this section we prove the two assertions in
Theorem~\ref{thm:confined-negative}. The proof is based on calculating the quantities $\Ea$ in \eqref{eq:E_functional} and $\Da$
in \eqref{eq:D_functional} along a suitably constructed sequence of densities $\{\rho_n\}_{n\in\N}$. The densities are constant multiples
of characteristic functions of sets with many small disjoint components. These sets are multi-dimensional analogs of the Dirac comb. Both estimates use the functional $J_\alpha$ defined in~\eqref{eq:J_functional}.

We proceed in three steps: (i) introduce the sequence $\{\rho_n\}_{n\in\N}$; (ii) establish auxiliary estimates common to the calculation of both functionals; (iii) complete the calculations
and obtain counterexamples to the lower bounds \eqref{eq:introlbineq} in Lebesgue classes with $p,q>1$. The determinant counterexample requires one additional estimate, which allows us to include the endpoint $q=1$ when $d\ge2$.

\subsection{Structure of the example}\label{seq:exampleden}
Fix $d\in\N$ and $1\le q<\infty$. For $d\ge2$, set
\begin{equation}\label{eq:set_V_d}
V_d:=[1/4,3/4]^{d-1}.
\end{equation}
Let $\{a_n\}_{n\in\N}$ be a sequence of amplitudes satisfying the
condition below, and define $\delta_n$ and $\ell_n$ as follows:
\begin{equation}\label{eq:layered_parameters}
1<a_n<a_{n+1},
\qquad
\delta_n:=a_n^{-q},
\qquad
\ell_n:=\frac{\delta_n}{n}.
\end{equation}
It is clear that
\[
a_n>1,
\qquad
0<\delta_n<1,
\qquad
0<\ell_n<\frac1n.
\]
For $k=0,\dots,n-1$, set
\begin{equation}\label{eq:layered_components}
s_{n,k}
:=
\frac1n\bigg(k+\frac12\bigg),
\qquad
I_{n,k}
:=
\bigg[ s_{n,k} -\frac{\ell_n}{10}, \, s_{n,k} + \frac{\ell_n}{10}\bigg],
\end{equation}
and define
\begin{equation}\label{eq:layered_sets}
A_n
:=
\bigcup_{k=0}^{n-1}I_{n,k},
\qquad
U_{n,d}
:=
\begin{dcases}
A_n,
\quad &d=1,\\[1mm]
A_n\times V_d,
\quad &d\ge2.
\end{dcases}
\end{equation}
The centerpoints $s_{n,k}$ are separated by $1/n$ and each interval $I_{n,k}$ has
half-length $\ell_n/10$. Thus the intervals $I_{n,0},\dots,I_{n,n-1}$ are pairwise
disjoint. Moreover,
\[
\begin{aligned}
\inf I_{n,0}
=
\frac1{2n}-\frac{\ell_n}{10}
&\ge
\frac1{2n}-\frac1{10n}
>0,
\\
\sup I_{n,n-1}
=
1-\frac1{2n}+\frac{\ell_n}{10}
&\le
1-\frac1{2n}+\frac1{10n}
<1.
\end{aligned}
\]
Hence
\[
A_n\subseteq(0,1),\qquad U_{n,d} \subseteq (0,1)^d.
\]

\begin{figure}[H]
\centering

\begin{minipage}[b]{0.40\textwidth}
\centering
\begin{tikzpicture}[scale=4.30]

    \def\nblocks{6}

    \def\w{0.040}

    \foreach \k in {0,...,5}{
        \pgfmathsetmacro{\x}{(\k+0.5)/\nblocks}

        \fill[gray!65]
            (\x-\w,0.25) rectangle (\x+\w,0.75);

        \draw[gray!72]
            (\x-\w,0.25) rectangle (\x+\w,0.75);
    }

    \draw[dashed,thick,gray!55]
        (0,0) rectangle (1,1);

    \node at (0.50,0.88)
        {\scriptsize $U_{n,2}=A_n\times V_2$};

    \node[anchor=west] at (1.00,0.50)
        {\scriptsize $V_2=\left[\tfrac14,\tfrac34\right]$};

\end{tikzpicture}

\par\smallskip
{\small (a) The set $U_{n,2}\subseteq(0,1)^2$.}
\end{minipage}
\hspace{0\textwidth}
\begin{minipage}[b]{0.56\textwidth}
\centering
\begin{tikzpicture}[
    x={(1.15cm,0cm)},
    y={(0.42cm,0.26cm)},
    z={(0cm,1.15cm)},
    scale=3.05
]

    \def\nblocks{6}

    \def\w{0.025}

    \draw[dashed,thick,gray!55] (0,1,0) -- (0,1,1); 
    \draw[dashed,thick,gray!55] (0,1,1) -- (1,1,1); 
    \draw[dashed,thick,gray!55] (0,1,0) -- (1,1,0); 
    \draw[dashed,thick,gray!55] (0,0,1) -- (0,1,1); 
    \draw[dashed,thick,gray!55] (1,1,0) -- (1,1,1); 

    \foreach \k in {0,...,5}{
        \pgfmathsetmacro{\x}{(\k+0.5)/\nblocks}
        \pgfmathsetmacro{\xa}{\x-\w}
        \pgfmathsetmacro{\xb}{\x+\w}

        \fill[gray!58]
            (\xa,0.25,0.25) --
            (\xb,0.25,0.25) --
            (\xb,0.25,0.75) --
            (\xa,0.25,0.75) -- cycle;

        \fill[gray!35]
            (\xa,0.25,0.75) --
            (\xb,0.25,0.75) --
            (\xb,0.75,0.75) --
            (\xa,0.75,0.75) -- cycle;

        \fill[gray!72]
            (\xb,0.25,0.25) --
            (\xb,0.75,0.25) --
            (\xb,0.75,0.75) --
            (\xb,0.25,0.75) -- cycle;

        \draw[gray!72]
            (\xa,0.25,0.25) --
            (\xb,0.25,0.25) --
            (\xb,0.25,0.75) --
            (\xa,0.25,0.75) -- cycle;

        \draw[gray!72]
            (\xa,0.25,0.75) --
            (\xb,0.25,0.75) --
            (\xb,0.75,0.75) --
            (\xa,0.75,0.75) -- cycle;

        \draw[gray!72]
            (\xb,0.25,0.25) --
            (\xb,0.75,0.25) --
            (\xb,0.75,0.75) --
            (\xb,0.25,0.75) -- cycle;
    }


\draw[dashed,thick,gray!55]
    (0,0,0) -- (1,0,0); 

\draw[dashed,thick,gray!55]
    (0,0,0) -- (0,0,1); 

\draw[dashed,thick,gray!55]
    (0,0,1) -- (1,0,1); 

\draw[dashed,thick,gray!55]
    (1,0,0) -- (1,0,1); 

\draw[dashed,thick,gray!55]
    (0,0,0) -- (0,1,0); 

\draw[dashed,thick,gray!55]
    (1,0,0) -- (1,1,0); 

\draw[dashed,thick,gray!55]
    (1,0,1) -- (1,1,1); 
    \node at (0.52,0.50,0.965)
        {\scriptsize $U_{n,3}=A_n\times V_3$};

    \node[anchor=east] at (-0.16,0.50,0.50)
        {\scriptsize $V_3=\left[\tfrac14,\tfrac34\right]^2$};

\end{tikzpicture}

\par\smallskip
{\small (b) The set $U_{n,3}\subseteq(0,1)^3$.}
\end{minipage}

\caption{Schematic illustration of the sets $U_{n,d}$ for $d=2$ and
$d=3$. In both panels, the dark regions represent the components
$I_{n,k}\times V_d$, where the intervals $I_{n,k}$ are centered at the
points $s_{n,k}=(k+\tfrac12)/n$. The drawings are not to scale: the true
half-length $\ell_n/10$ is enlarged for visibility.}
\label{fig:layered-sets}
\end{figure}
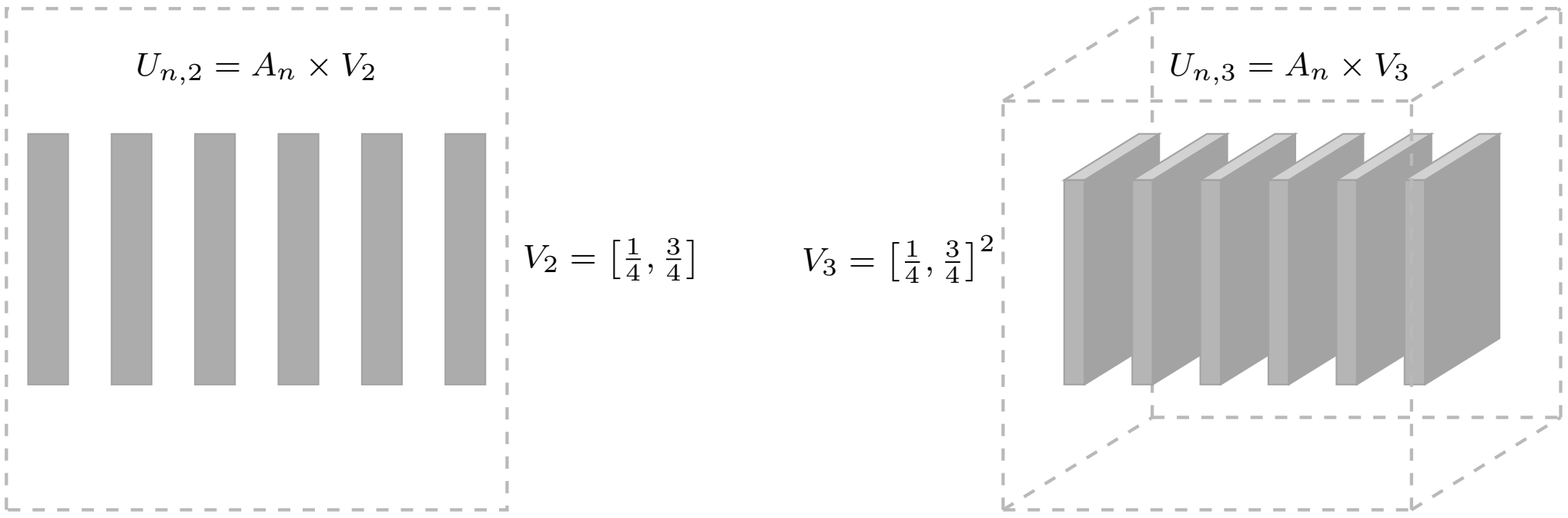

This procedure defines the sets $U_{n,d}$; let $\chi_{U_{n,d}}$ denote the associated characteristic function.
The densities $\rho_n$ are defined by
\begin{equation}
\rho_n := a_n \chi_{U_{n,d}},
\end{equation}
where the amplitudes $a_n$ will be appropriately selected at the end.
By construction $\supp \rho_n \subseteq (0,1)^d$ and, for the chosen exponent $1 \leq q < \infty$, we have
\begin{equation}
\| \rho_n \|_{L^q(\R^d)}^q = a_n^q \, |U_{n,d}| = c_d,
\end{equation}
where $c_d$ is independent of $n$ and the amplitude $a_n$.  The next objective is to calculate
$J_\alpha ( \chi_{U_{n,d}})$, where $J_\alpha$ is the functional defined in~\eqref{eq:J_functional}. This will be based on the auxiliary calculations
presented below.

\subsection{An auxiliary kernel}
For
$0<\alpha<d$ and $r\neq0$, define
\begin{equation}\label{eq:H_alpha,d}
H_{\alpha,d}(r)
:=
\begin{dcases}
|r|^{\alpha-1},
\qquad &d=1,\\[2mm]
\displaystyle
\int_{V_d}\int_{V_d}
\big(r^2+|u-v|^2\big)^{\frac{\alpha-d}{2}}
\,\du\,\dv,
\qquad &d\ge2.
\end{dcases}
\end{equation}
The first lemma provides bounds on the function $H_{\alpha,d}(r)$ that will
be the cornerstone of forthcoming estimations.

\begin{lemma}
\label{lem:transverse-kernel-estimate}
Let $d\in\N$ and $0<\alpha<d$. Then, for $0<|r|\le1$,
\begin{equation}\label{eq:H_alpha,d_bounds}
H_{\alpha,d}(r)
\le
\begin{dcases}
c_{\alpha,d} \, |r|^{\alpha-1},
\qquad &0<\alpha<1,\\[1mm]
c_d\big(1+\big|\!\log|r|\big|\big),
\qquad &\alpha=1,\\[1mm]
c_{\alpha,d},
\qquad &1<\alpha<d.
\end{dcases}
\end{equation}
\end{lemma}

\begin{proof}
If $d=1$, then $0<\alpha<1$ and $H_{\alpha,1}(r)=|r|^{\alpha-1}$. 
Thus the first estimate holds with constant one.

Assume that $d\ge2$. Since $H_{\alpha,d}$ is even, it is enough to consider $0<r\le1$. Changing variables from $(u,v)$ to $(u,\eta)$, where $\eta=v-u$, and using Tonelli's theorem, we have
\begin{align*}
H_{\alpha,d}(r)
&=
\int_{V_d}
\int_{V_d-u}
\big(r^2+|\eta|^2\big)^{\frac{\alpha-d}{2}}
\,\rd\eta\,\du
\\
&\leq 
\int_{V_d}
\int_{V_d-V_d}
\big(r^2+|\eta|^2\big)^{\frac{\alpha-d}{2}}
\,\rd\eta\,\du
\\
&=|V_d|
\int_{V_d-V_d}
\big(r^2+|\eta|^2\big)^{\frac{\alpha-d}{2}}\,\rd\eta.
\end{align*}
Since $V_d-V_d=[-1/2,1/2]^{d-1}$,  we have
\[
V_d-V_d
\subseteq
\overline{B_{\kappa_d}^{d-1}},
\qquad
\kappa_d
:=
\max\bigg\{
1,\frac{\sqrt{d-1}}2
\bigg\}.
\]
Using polar coordinates in $\R^{d-1}$, we obtain
\[
H_{\alpha,d}(r)
\le
c_d
\int_0^{\kappa_d}
\big(r^2+s^2\big)^{\frac{\alpha-d}{2}}
s^{d-2}\,\ds.
\]

It remains to calculate the integral on the right-hand side.
Using $\kappa_d\ge1$ and $0<r\le1$, we split the integral
\begin{align*}
\int_0^{\kappa_d}
\big(r^2+s^2\big)^{\frac{\alpha-d}{2}}
s^{d-2}\,\ds
&=
\int_0^r
\big(r^2+s^2\big)^{\frac{\alpha-d}{2}}
s^{d-2}\,\ds
+
\int_r^{\kappa_d}
\big(r^2+s^2\big)^{\frac{\alpha-d}{2}}
s^{d-2}\,\ds.
\end{align*}
For the first integral, using $\big(r^2+s^2\big)^{\frac{\alpha-d}{2}} \le r^{\alpha-d}$,
we deduce
\[
\int_0^r
\big(r^2+s^2\big)^{\frac{\alpha-d}{2}}
s^{d-2}\,\ds
\le
r^{\alpha-d}
\int_0^r s^{d-2}\,\ds
=
c_d \, r^{\alpha-1}.
\]
For the second integral, using $\big(r^2+s^2\big)^{\frac{\alpha-d}{2}} \le s^{\alpha-d}$, we obtain
\[
\int_r^{\kappa_d}
\big(r^2+s^2\big)^{\frac{\alpha-d}{2}}
s^{d-2}\,\ds
\le
\int_r^{\kappa_d}s^{\alpha-2}\,\ds.
\]
To conclude, note that if $0<\alpha<1$, then
\[
\int_r^{\kappa_d}s^{\alpha-2}\,\ds
=
\frac{r^{\alpha-1}-\kappa_d^{\alpha-1}}{1-\alpha}
\le
c_\alpha \, r^{\alpha-1}.
\]
If $\alpha=1$, then
\[
\int_r^{\kappa_d}\frac{1}{s} \, \ds
=
\log\frac{\kappa_d}{r}
\le
c_d\big(1+|\!\log r|\big).
\]
Finally, if $1<\alpha<d$, then
\[
\int_r^{\kappa_d}s^{\alpha-2}\,\ds
=
\frac{\kappa_d^{\alpha-1}-r^{\alpha-1}}{\alpha-1}
\le
c_{\alpha,d}.
\]

Combining the estimates above proves the stated bounds.
\end{proof}

We next express $J_\alpha(\chi_{U_{n,d}})$ where $J_\alpha$ is defined in~\eqref{eq:J_functional}. For $d\geq2$, setting $z = (t,u)$, $w=(s,v)$, $t,s \in A_n$, $u, v \in V_d$,  a direct calculation using the definition of $U_{n,d}$ and $H_{\alpha,d}$ shows
\begin{equation}\label{eq:layered-pair-reduction}
\begin{split}
J_\alpha(\chi_{U_{n,d}})
&=
\int_{U_{n,d}} \int_{U_{n,d}} \Big ( |t-s|^2 + |u-v|^2 \Big )^\frac{\alpha-d}{2} \chi_{A_n}(t) \chi_{V_d}(u) 
\chi_{A_n}(s) \chi_{V_d}(v) \, \dz \, \dw
\\
&= \int_{A_n} \int_{A_n} \int_{V_d} \int_{V_d} \Big ( |t-s|^2 + |u-v|^2 \Big )^\frac{\alpha-d}{2} \, \du \, \dv \, \dt \, \ds
\\
& =
\int_{A_n}\int_{A_n}
H_{\alpha,d}(s-t)\,\dt\,\ds \\
&= \sum_{0\le k,m\le n-1}
\int_{I_{n,k}}\int_{I_{n,m}}
H_{\alpha,d}(s-t)\,\dt\,\ds.
\end{split}
\end{equation}
For $d=1$, the last two expressions follow directly from $U_{n,1}=A_n$ and $H_{\alpha,1}(r)=|r|^{\alpha-1}$.

Since the intervals $I_{n,k}$ are pairwise disjoint, the last integral in~\eqref{eq:layered-pair-reduction} is split according to
whether~$s$ and~$t$ belong to the same interval or to different intervals. For $n\ge2$, define
\begin{align}
\label{eq:layered-self-interaction-definition}
J_\alpha^{\rm self}(U_{n,d})
&:=
\sum_{k=0}^{n-1}
\int_{I_{n,k}}\int_{I_{n,k}}
H_{\alpha,d}(s-t)\,\dt\,\ds,
\\
\label{eq:layered-cross-interaction-definition}
J_\alpha^{\rm cross}(U_{n,d})
&:=
\sum_{\substack{0\le k,m\le n-1\\k\neq m}}
\int_{I_{n,k}}\int_{I_{n,m}}
H_{\alpha,d}(s-t)\,\dt\,\ds,
\end{align}
and use these to express
\begin{equation}\label{eq:splitting}
J_\alpha(\chi_{U_{n,d}}) = J_\alpha^{\rm self}(U_{n,d}) +  J_\alpha^{\rm cross}(U_{n,d}).
\end{equation}
The next lemma provides the required estimates for the self-interacting and
cross-interacting terms.

\begin{lemma}
\label{lem:layered-pair-interaction}
Let $d\in\N$ and $0<\alpha<d$. Then, for every $n\ge2$, 
\begin{equation}\label{eq:layered-self-estimate}
J_\alpha^{\rm self}(U_{n,d})
\le
\begin{dcases}
c_{\alpha,d} \,
n^{-\alpha}\delta_n^{\alpha+1},
\qquad &0<\alpha<1,\\[1mm]
c_d \,
n^{-1}\delta_n^2
\bigg(
1+\log\frac{n}{\delta_n}
\bigg),
\qquad &\alpha=1,\\[2mm]
c_{\alpha,d} \,
n^{-1}\delta_n^2,
\qquad &1<\alpha<d,
\end{dcases}
\end{equation}
and
\begin{equation}\label{eq:layered-pair-estimate}
J_\alpha^{\rm cross}(U_{n,d})
\le
c_{\alpha,d} \, \delta_n^2.
\end{equation}
\end{lemma}

\begin{proof}
We first estimate the self-interaction terms. Set
$
L:=|I_{n,k}|=\frac{\ell_n}{5}.
$
Suppose that $0<\alpha<1$. By
Lemma~\ref{lem:transverse-kernel-estimate},
\[
H_{\alpha,d}(s-t)
\le
c_{\alpha,d} \, |s-t|^{\alpha-1}.
\]
Translating $I_{n,k}$ to $[0,L]$ and using symmetry, we obtain
\begin{align*}
\int_{I_{n,k}}\int_{I_{n,k}}
|s-t|^{\alpha-1}\,\dt\,\ds
&=
2\int_0^L\int_0^s
(s-t)^{\alpha-1}\,\dt\,\ds
\\
&=
\frac2\alpha\int_0^L s^\alpha\,\ds
\\
&=
\frac{2}{\alpha(\alpha+1)}
L^{\alpha+1}.
\end{align*}
Consequently,
\[
J_\alpha^{\rm self}(U_{n,d})
\le
c_{\alpha,d} \, n \, \ell_n^{\alpha+1}
=
c_{\alpha,d} \, n^{-\alpha}\delta_n^{\alpha+1}.
\]

Suppose next that $\alpha=1$. Then
\[
H_{1,d}(s-t)
\le
c_d\big(1+\big|\!\log|s-t|\big|\big).
\]
Since
$ L=\frac{\ell_n}{5}\le\frac1{10}$, 
we have
$
\big|\!\log|s-t|\big|=-\log|s-t|
$
whenever $s,t\in I_{n,k}$ and $s\neq t$. Translating again to $[0,L]$,
\begin{align*}
\int_{I_{n,k}}\int_{I_{n,k}}
\big(1+\big|\!\log|s-t|\big|\big)
\,\dt\,\ds &=
2\int_0^L\int_0^s
\big(1-\log(s-t)\big)
\,\dt\,\ds
\\
&=
2\int_0^L
\big(2s-s\log s\big)\,
\ds
\\
&=
L^2\bigg(\frac52-\log L\bigg)
\\
&=
\bigg(\frac{\ell_n}{5}\bigg)^2
\bigg(
\frac52+\log\frac5{\ell_n}
\bigg).
\end{align*}
It follows that
\begin{align*}
J_1^{\rm self}(U_{n,d})
\le
c_d \, n \, \ell_n^2
\big(1+|\!\log\ell_n|\big)
=
c_d \, n^{-1}\delta_n^2
\bigg(
1+\log\frac n{\delta_n}
\bigg).
\end{align*}

Finally, if $1<\alpha<d$, then
\[
H_{\alpha,d}(s-t)\le c_{\alpha,d}.
\]
Therefore
\[
J_\alpha^{\rm self}(U_{n,d})
\le
c_{\alpha,d} \, n \, \ell_n^2
=
c_{\alpha,d} \, n^{-1}\delta_n^2.
\]
This proves~\eqref{eq:layered-self-estimate}.

Next, we estimate the cross-interaction terms. Let
$s\in I_{n,k}$,  $t\in I_{n,m}$,  $k\neq m$ and set
\[
j:=|k-m|\ge1.
\]
The distance between the centers of $I_{n,k}$ and $I_{n,m}$ is $j/n$. Since $s$ and $t$
differ from their respective centers by at most $\ell_n/10$,
\[
|s-t|
\ge
\frac jn-\frac{\ell_n}{5}.
\]
Using $\ell_n\le1/n$, we obtain
\begin{equation}\label{eq:layered-component-separation}
|s-t|
\ge
\frac1n\bigg(j-\frac15\bigg)
\ge
\frac{4j}{5n}.
\end{equation}

Suppose first that $0<\alpha<1$. By
Lemma~\ref{lem:transverse-kernel-estimate} and
\eqref{eq:layered-component-separation},
\[
H_{\alpha,d}(s-t)
\le
c_{\alpha,d}\,
n^{1-\alpha}|k-m|^{\alpha-1}.
\]
Using that $|I_{n,k}|=|I_{n,m}|=\frac{\ell_n}{5}$, 
symmetry gives
\begin{align*}
J_\alpha^{\rm cross}(U_{n,d})
&\le
c_{\alpha,d}\,
n^{1-\alpha}\ell_n^2
\sum_{0\le k<m\le n-1}
(m-k)^{\alpha-1}
\\
&=
c_{\alpha,d}\,
n^{1-\alpha}\ell_n^2
\sum_{k=0}^{n-2}
\sum_{m=k+1}^{n-1}
(m-k)^{\alpha-1}
\\
&=
c_{\alpha,d}\,
n^{1-\alpha}\ell_n^2
\sum_{k=0}^{n-2}
\sum_{j=1}^{n-1-k}
j^{\alpha-1}
\\
&=
c_{\alpha,d}\,
n^{1-\alpha}\ell_n^2
\sum_{j=1}^{n-1}
\sum_{k=0}^{n-1-j}
j^{\alpha-1}
\\
&=
c_{\alpha,d}\,
n^{1-\alpha}\ell_n^2
\sum_{j=1}^{n-1}
(n-j)j^{\alpha-1}.
\end{align*}
Moreover,
\[
\sum_{j=1}^{n-1}
(n-j)j^{\alpha-1}
\le
n\sum_{j=1}^{n-1}j^{\alpha-1}
\le
c_\alpha \, n^{\alpha+1},
\]
where we used
\[
j^{\alpha-1} = \int_{j-1}^j j^{\alpha-1}\,\dx
\le
\int_{j-1}^j x^{\alpha-1}\,\dx.
\]
Hence
\[
J_\alpha^{\rm cross}(U_{n,d})
\le
c_{\alpha,d} \, n^2\ell_n^2
=
c_{\alpha,d} \, \delta_n^2.
\]

Suppose next that $\alpha=1$. Since $A_n\subseteq(0,1)$, we have
$
0<|s-t|<1
$
when $s$ and $t$ belong to different components. Thus
\[
\big|\!\log|s-t|\big|
=
\log\frac1{|s-t|}.
\]
Using~\eqref{eq:layered-component-separation},
\[
\log\frac1{|s-t|}
\le
\log\frac{5n}{4j}
=
\log\frac54+\log\frac nj.
\]
Lemma~\ref{lem:transverse-kernel-estimate} therefore gives
\[
H_{1,d}(s-t)
\le
c_d
\bigg(
1+\log\frac nj
\bigg).
\]
By symmetry and the same reindexing used in the preceding case, with $j=m-k$, we obtain
\[
J_1^{\rm cross}(U_{n,d})
\le
c_d \, \ell_n^2
\sum_{j=1}^{n-1}
(n-j)
\bigg(
1+\log\frac nj
\bigg).
\]
The nonlogarithmic part is bounded by $n^2$. For the logarithmic part,
\[
\sum_{j=1}^{n-1}
(n-j)\log\frac nj
\le
n\sum_{j=1}^{n-1}\log\frac nj.
\]
Since the function $x\mapsto\log(n/x)$ is decreasing, for every
$j=2,\dots,n-1$ we have
\[
\log\frac nj
=
\int_{j-1}^j\log\frac nj\,\dx
\le
\int_{j-1}^j\log\frac nx\,\dx.
\]
Therefore
\begin{align*}
\sum_{j=1}^{n-1}\log\frac nj
&=
\log n+\sum_{j=2}^{n-1}\log\frac nj
\\
&\le
\log n+\int_1^n\log\frac nx\,\dx
\\
&=
n-1.
\end{align*}
It follows that
\[
J_1^{\rm cross}(U_{n,d})
\le
c_d \, n^2\ell_n^2
=
c_d \, \delta_n^2.
\]

Finally, suppose that $1<\alpha<d$. By
Lemma~\ref{lem:transverse-kernel-estimate},
\[
H_{\alpha,d}(s-t)\le c_{\alpha,d}.
\]
Since the intervals $I_{n,k}$ are pairwise disjoint,
\[
|A_n|
=
\sum_{k=0}^{n-1}|I_{n,k}|
=
\frac{n\ell_n}{5}.
\]
Therefore
\[
J_\alpha^{\rm cross}(U_{n,d})
\le
c_{\alpha,d} \, |A_n|^2
\le
c_{\alpha,d} \,(n\ell_n)^2
=
c_{\alpha,d} \, \delta_n^2.
\]
This proves~\eqref{eq:layered-pair-estimate}. Combining
\eqref{eq:layered-self-estimate} and
\eqref{eq:layered-pair-estimate} gives
\eqref{eq:layered-full-pair-estimate}.
\end{proof}

In summary, using Lemma \ref{lem:layered-pair-interaction}, we obtain an estimate of the strength
of interactions in \eqref{eq:layered-pair-reduction} in terms of the parameter $\delta_n$:
\begin{equation}\label{eq:layered-full-pair-estimate}
J_\alpha(\chi_{U_{n,d}})
\le
\begin{dcases}
c_{\alpha,d}
\left(
n^{-\alpha}\delta_n^{\alpha+1}
+
\delta_n^2
\right),
\qquad &0<\alpha<1,\\[2mm]
c_d
\bigg(
n^{-1}\delta_n^2
\Big(
1+\log\frac{n}{\delta_n}
\Big)
+
\delta_n^2
\bigg),
\qquad &\alpha=1,\\[2mm]
c_{\alpha,d}
\left(
n^{-1}\delta_n^2
+
\delta_n^2
\right),
\qquad &1<\alpha<d.
\end{dcases}
\end{equation}

\subsection{Proof of Theorem~\ref{thm:confined-negative} part (i)}
\begin{proof}
Let $\rho_n = a_n \, \chi_{U_{n,d}}$, where $a_n$ and $U_{n,d}$ are defined in~\eqref{eq:layered_parameters}
and~\eqref{eq:layered_sets}, respectively.
Recall that  $\supp \rho_n \subseteq (0,1)^d$ and that $\|\rho_n\|_{L^q} = c_{d,q}$. We still have the freedom to select $a_n$ subject to the conditions that $\{ a_n \}_{n \in \N}$ is increasing and $a_n > 1$.
It remains to prove the asserted convergence of $\Ea(\rho_n)$.

We have
\begin{equation}\label{eq:energy-layered-domination}
J_\alpha(\rho_n)
=
a_n^2 \, J_\alpha(\chi_{U_{n,d}}),
\end{equation}
where $J_\alpha$ is defined in~\eqref{eq:J_functional}.
Using~\eqref{eq:layered-full-pair-estimate} in
\eqref{eq:energy-layered-domination}, and recalling from
\eqref{eq:layered_parameters} that
$\delta_n=a_n^{-q}$, 
we obtain
\[
J_\alpha(\rho_n)
\le
\begin{dcases}
c_{\alpha,d}
\big(
n^{-\alpha}a_n^{2-q(\alpha+1)}
+
a_n^{2-2q}
\big),
\qquad&0<\alpha<1,\\[2mm]
c_d
\Big(
n^{-1}a_n^{2-2q}
\big(1+\log(n \,a_n^q)\big)
+
a_n^{2-2q}
\Big),
\qquad&\alpha=1,\\[2mm]
c_{\alpha,d}
\big(
n^{-1}a_n^{2-2q}
+
a_n^{2-2q}
\big),
\qquad&1<\alpha<d.
\end{dcases}
\]

Since $q>1$, for any  $a_n\to\infty$,
\[
a_n^{2-2q}\to0 \qquad \text{as } n \to \infty.
\]
Moreover, if we now select
\begin{equation}\label{basichyp}
a_n=\log(n+2),
\end{equation}
then
\[
n^{-\sigma}a_n^\gamma\to0
\qquad
\text{as }n\to\infty
\]
for every $\sigma>0$ and every $\gamma\in\R$. This proves the required convergence to zero of
the first term when $0<\alpha<1$ and when $1<\alpha<d$.

It remains to consider the first term when $\alpha=1$ for the selection \eqref{basichyp}. For $n\ge1$,
\[
1+\log(n \, a_n^q)
=
1+\log n+q\log a_n
\le
c_q \, a_n.
\]
Therefore
\[
n^{-1}a_n^{2-2q}
\big(1+\log(n \, a_n^q)\big)
\le
c_q \, n^{-1}a_n^{3-2q}
\to0 \qquad
\text{as }n\to\infty.
\]
Thus, in every case,
\[
J_\alpha(\rho_n)\to0
\qquad
\text{as }n\to\infty.
\]
Since $\mathcal E_\alpha(\rho_n) = \tfrac1{d-\alpha}J_\alpha(\rho_n)$, 
the result follows.
\end{proof}

\subsection{The determinant reduction}

Next, we focus on the functional $\Da$, again for densities of the form $0 \le \rho \le a \chi_U$. In preparation, we show that estimating the determinant functional reduces to controlling
\begin{equation}\label{eq:mathcal_I_alpha}
\mathcal I_\alpha(U)
:=
\int_0^1
\int_{\R^d}
\int_{\R^d}
\chi_U(z)\,\chi_U(w)\,
\chi_U \big(\theta z+(1-\theta)w\big)\,
|z-w|^{\alpha-d}
\,\dz\,\dw\,\rd\theta ,
\end{equation}
where $U\subseteq\R^d$ is a bounded measurable set.

\begin{lemma}
\label{lem:trace-incidence-bound}
Let $d\in\N$, $0<\alpha<d$, and $\rho$ be a density satisfying
\[
0\le\rho\le a \,\chi_U,
\]
for some $a>0$ and some bounded measurable set $U\subseteq\R^d$. Then
\[
\Da(\rho)
\le
c_d\,
a^{2+1/d} \, \mathcal I_\alpha(U).
\]
\end{lemma}

\begin{proof}
Applying the arithmetic-geometric mean inequality to the eigenvalues of $\Sa(\rho)(x)$, we obtain
\[
\Da(\rho)
\le
\frac1d
\int_{\R^d}
\rho(x)^{1/d}\tr\Sa(\rho)(x)\,\dx.
\]
Taking the trace of $\Sa(\rho)(x)$ and using
$\rho\leq a\chi_U$ gives
\begin{align*}
\tr\Sa(\rho)(x)
&=
\frac12
\int_0^1
\int_{\R^d}
\rho\big(x+(\theta-1)y\big)\,
\rho\big(x+\theta y\big)\,
|y|^{\alpha-d}
\,\dy\,\rd\theta
\\
&\le
\frac{a^2}{2}
\int_0^1
\int_{\R^d}
\chi_U\big(x+(\theta-1)y\big)\,
\chi_U\big(x+\theta y\big)\,
|y|^{\alpha-d}
\,\dy\,\rd\theta.
\end{align*}
Therefore
\begin{align*}
\Da(\rho)
&\le
\frac{1}{2d}\,a^{2+1/d}
\int_0^1
\int_{\R^d}
\int_{\R^d}
\chi_U(x)\,
\chi_U\big(x+(\theta-1)y\big) \,
\chi_U\big(x+\theta y\big)\,
|y|^{\alpha-d}
\,\dy\,\dx\,\rd\theta.
\end{align*}
For fixed $\theta$, the change of variables
\[
z=x+(\theta-1)y,
\qquad
w=x+\theta y,
\]
has inverse
\[
y=w-z,
\qquad
x=\theta z+(1-\theta)w,
\]
and Jacobian one; it yields the desired result.
\end{proof}

The critical step consists of using the self- and cross-interaction
estimates in Lemma~\ref{lem:layered-pair-interaction} to estimate
$\mathcal I_\alpha(U_{n,d})$.

\begin{lemma}
\label{lem:layered-three-point-interaction}
Let $d \in \N$ and $0<\alpha<d$. Then, for every integer $n\ge2$,
\begin{equation}\label{eq:estimationI}
\mathcal I_\alpha(U_{n,d})
\le
\begin{dcases}
c_{\alpha,d}
\big(
n^{-\alpha}\delta_n^{\alpha+1}
+
\delta_n^3
\big),
\qquad &0<\alpha<1,\\[2mm]
c_d
\bigg(
n^{-1}\delta_n^2
\Big(
1+\log\frac{n}{\delta_n}
\Big)
+
\delta_n^3
\bigg),
\qquad &\alpha=1,\\[2mm]
c_{\alpha,d}
\big(
n^{-1}\delta_n^2
+
\delta_n^3
\big),
\qquad &1<\alpha<d.
\end{dcases}
\end{equation}
\end{lemma}

\begin{proof}
Suppose first that $d\ge2$, and write
\[
z=(s,u),
\qquad
w=(t,v), \qquad s,t \in \R, \qquad u,v \in \R^{d-1}.
\]
Since
\[
U_{n,d}=A_n\times V_d
\]
and $V_d$ is convex, we have
\[
\theta u+(1-\theta)v\in V_d
\qquad
\text{for all }u,v\in V_d
\text{ and }\theta\in[0,1].
\]
Tonelli's theorem and the definition of $H_{\alpha,d}$ therefore give
\begin{equation}\label{eq:layered-three-point-reduction}
\mathcal I_\alpha(U_{n,d})
=
\int_{A_n}\int_{A_n}
H_{\alpha,d}(s-t)
\left(
\int_0^1
\chi_{A_n}\big(\theta s+(1-\theta)t\big)
\,\rd\theta
\right)
\dt\,\ds.
\end{equation}
When $d=1$, the same identity follows directly from
\[
H_{\alpha,1}(s-t)=|s-t|^{\alpha-1}.
\]
The bounds in Lemma~\ref{lem:transverse-kernel-estimate} imply that
$H_{\alpha,d}$ is locally integrable near the origin, so the right-hand side is integrable.
Also, in analogy to \eqref{eq:splitting}, $\mathcal I_\alpha(U_{n,d})$ can be split into 
contributions of self- and cross-interacting terms.

Suppose that $s,t\in I_{n,k}$ for some $k\in\{0,\ldots,n-1\}$. Since $I_{n,k}$ is an interval,
\[
\theta s+(1-\theta)t\in I_{n,k}\subseteq A_n
\qquad
\text{for all }\theta\in[0,1].
\]
Hence
\[
\int_0^1
\chi_{A_n}\big(\theta s+(1-\theta)t\big)
\,\rd\theta
=
1.
\]
The part of~\eqref{eq:layered-three-point-reduction} corresponding to pairs in the same
interval is therefore equal to
\[
J_\alpha^{\rm self}(U_{n,d}),
\]
and is bounded by~\eqref{eq:layered-self-estimate}.

We now consider pairs in different intervals. Let
$s\in I_{n,k}$, $t\in I_{n,m}$ with $k\neq m$  and set $j:=|k-m|\ge1$. 
By~\eqref{eq:layered-component-separation},
\[
|s-t|
\ge
\frac{4j}{5n}.
\]
Since $s\neq t$, the change of variables
\[
x=\theta s+(1-\theta)t
=
t+\theta(s-t)
\]
gives
\begin{align*}
\int_0^1
\chi_{A_n}\big(\theta s+(1-\theta)t\big)
\,\rd\theta
&=
\frac1{|s-t|}
\int_{\min\{s,t\}}^{\max\{s,t\}}
\chi_{A_n}(x)\,\dx
\\
&=
\frac{1}{|s-t|}\big|A_n\cap[\min\{s,t\},\max\{s,t\}]\big|.
\end{align*}
Moreover,
\[
A_n\cap[\min\{s,t\},\max\{s,t\}]
\subseteq
\bigcup_{r=\min\{k,m\}}^{\max\{k,m\}}I_{n,r}.
\]
The union on the right consists of $j+1$ intervals, each of length $\ell_n/5$. Hence
\[
\big|A_n\cap[\min\{s,t\},\max\{s,t\}]\big|
\le
\frac{(j+1)\ell_n}{5}.
\]
Using $n\,\ell_n=\delta_n$, we obtain
\begin{align*}
\int_0^1
\chi_{A_n}\big(\theta s+(1-\theta)t\big)
\,\rd\theta
&\le
\frac{(j+1)\ell_n/5}{4j/(5n)}
=
\frac{j+1}{4j}\,\delta_n
\le
\frac{\delta_n}{2}.
\end{align*}
Consequently, the part of~\eqref{eq:layered-three-point-reduction} corresponding to pairs in
different intervals is bounded by $\frac{\delta_n}{2}J_\alpha^{\rm cross}(U_{n,d})$. 
By~\eqref{eq:layered-pair-estimate},
\begin{equation}\label{eq:cross-interactions}
\frac{\delta_n}{2}
J_\alpha^{\rm cross}(U_{n,d})
\le
c_{\alpha,d} \,\delta_n^3 .
\end{equation}
Combining \eqref{eq:cross-interactions} and \eqref{eq:layered-self-estimate} proves \eqref{eq:estimationI}.
\end{proof}

\subsection{Proof of Theorem~\ref{thm:confined-negative} part (ii)}
\begin{proof}
Let $\rho_n = a_n \, \chi_{U_{n,d}}$, where $a_n$ and $U_{n,d}$ are defined in~\eqref{eq:layered_parameters}
and~\eqref{eq:layered_sets}, respectively. The support and norm conditions follow from the construction.
It remains to prove the asserted convergence of $\Da(\rho_n)$. We again select $a_n$ to satisfy \eqref{basichyp}. 

Lemma~\ref{lem:trace-incidence-bound} implies
\[
\Da(\rho_n)
\le
c_d\,a_n^\beta \, \mathcal I_\alpha(U_{n,d}).
\]
Applying Lemma~\ref{lem:layered-three-point-interaction} and using
$\delta_n=a_n^{-q}$, we obtain
\begin{equation}\label{eq:confined-obstruction-final-estimates}
\Da(\rho_n)
\le
\begin{dcases}
c_{\alpha,d}
\big(
n^{-\alpha}a_n^{\beta-q(\alpha+1)}
+
a_n^{\beta-3q}
\big),
&\qquad 0<\alpha<1,\\[2mm]
c_d
\bigg(
n^{-1}a_n^{\beta-2q}
\Big(1+\log(n \, a_n^q)\Big)
+
a_n^{\beta-3q}
\bigg),
&\qquad \alpha=1,\\[2mm]
c_{\alpha,d}
\big(
n^{-1}a_n^{\beta-2q}
+
a_n^{\beta-3q}
\big),
&\qquad 1<\alpha<d.
\end{dcases}
\end{equation}
We first consider the second term in each line of \eqref{eq:confined-obstruction-final-estimates}. Recall that
$\beta = 2 + \frac{1}{d}$. If $d\ge2$ and $q\ge1$, then
$\beta-3q \le -1+\frac1d <0$. 
If $d=1$ and $q>1$ then $\beta-3q =3(1-q)<0$.
For $a_n\to\infty$, in either case
\begin{equation}\label{eq:determinant-common-term-decay}
a_n^{\beta-3q}\to0
\qquad
\text{as }n\to\infty.
\end{equation}

We next consider the first term in each line of
\eqref{eq:confined-obstruction-final-estimates}. As in the proof of part~(i), we use the fact that
for $a_n = \log (n+2)$ we have 
$n^{-\sigma}a_n^\gamma\to0 \; \text{as }n\to\infty$
for every $\sigma>0$ and every $\gamma\in\R$. It follows that
\[
n^{-\alpha}a_n^{\beta-q(\alpha+1)}
\to0 \qquad \text{as } n \to \infty
\]
when $0<\alpha<1$, and that
\[
n^{-1}a_n^{\beta-2q}\to0 \qquad \text{as } n \to \infty
\]
when $1<\alpha<d$.
It remains to consider the first term when $\alpha=1$. As in the proof of part~(i),
again for $a_n = \log (n+2)$, we have
\[
1+\log(n \, a_n^q)
\le
c_q\,a_n
\qquad
\text{for all }n\ge1.
\]
Therefore
\[
n^{-1}a_n^{\beta-2q}
\big(1+\log(n \, a_n^q)\big)
\le
c_q\,n^{-1}a_n^{\beta-2q+1}
\to0 \qquad \text{as } n \to \infty.
\]
Combining these estimates with
\eqref{eq:determinant-common-term-decay} and
\eqref{eq:confined-obstruction-final-estimates} proves that $\Da(\rho_n) \to 0$.
\end{proof}

\section{Lower bounds in \texorpdfstring{$L^1$}{L1} under confinement}
\label{sec:lower-bounds-conf-alone}

The negative results for confined densities leave two endpoint cases to consider. For the energy functional $\Ea$, confinement alone yields an $L^1$ lower bound in every dimension. For the determinant functional $\Da$, the analogous endpoint bound holds in dimension one. We prove these two results below.

\begin{proposition}
\label{prop:E_L1_confined}
Let $d\in\N$, $0<\alpha<d$, and $R>0$. There exists a positive constant $C=C(\alpha,d,R)$ such that
\[
\mathcal E_\alpha(\rho)
\ge
C \,
\|\rho\|_{L^1(B_R)}^2
\]
for all densities $\rho$ supported in $B_R$.
\end{proposition}

\begin{proof}
If $\mathcal E_\alpha(\rho)=\infty$, there is nothing to prove. We therefore assume that
$\mathcal E_\alpha(\rho)<\infty$.

Since $\rho=0$ almost everywhere in $\R^d\setminus B_R$, the double integral may be restricted to $B_R\times B_R$. For $x,y\in B_R$, we have
\[
|x-y|\le2R,
\]
and so
\[
|x-y|^{\alpha-d}
\ge
(2R)^{\alpha-d},
\]
since $\alpha-d<0$. Therefore,
\begin{align*}
\mathcal E_\alpha(\rho)
&=
\frac1{d-\alpha}
\IR \IR
\rho(x) \, \rho(y) \, |x-y|^{\alpha-d}\,\dx\,\dy
\\
&\ge
\frac{(2R)^{\alpha-d}}{d-\alpha}
\int_{B_R} \int_{B_R}
\rho(x) \, \rho(y)\,\dx\,\dy
\\
&=
\frac{(2R)^{\alpha-d}}{d-\alpha}
\|\rho\|_{L^1(B_R)}^2.
\end{align*}
\end{proof}

\begin{proposition}
\label{prop:D-confined-L1-endpoint}
Let $d=1$, $0<\alpha<1$ and $R>0$. There exists a positive constant $C = C(\alpha,R)$ such that
\[
\Da(\rho)
\ge
C \,
\|\rho\|_{L^1(B_R)}^3
\]
for all densities $\rho$ supported in $B_R$.
\end{proposition}

\begin{proof}

If $\Da(\rho)=\infty$, there is nothing to prove. We therefore assume that $\Da(\rho)<\infty$. For $N\in\N$, set
\[
\rho_N:=\min\{\rho,N\}.
\]
Then $\rho_N$ is supported in $B_R$ and belongs to $L^1(\R) \cap L^\infty(\R)$. Since $d=1$,
the quantity $\Sa(\rho)$ is scalar. Moreover, $0\le\rho_N\le\rho$ and the definition of
$\Sa$ give
\[
0\le\Sa(\rho_N)(x)\le\Sa(\rho)(x),
\]
for almost every $x\in\R$. Consequently,
\begin{equation}\label{eq:one-dimensional-D-monotonicity}
\Da(\rho_N)
=
\int_{\R}\rho_N(x)\Sa(\rho_N)(x)\,\dx
\le
\int_{\R}\rho(x)\Sa(\rho)(x)\,\dx
=
\Da(\rho).
\end{equation}

We now estimate $\Da(\rho_N)$. Using the definition of $\Sa$ and the change of variables
\[
z=x+(\theta-1)y,
\qquad
w=x+\theta y,
\]
we obtain
\[
\Da(\rho_N)
=
\frac12
\int_0^1\int_{\R}\int_{\R}
\rho_N\big(\theta z+(1-\theta)w\big)\,
\rho_N(z) \, \rho_N(w) \, |z-w|^{\alpha-1}
\,\dz\,\dw\,\rd\theta.
\]
For $z\neq w$,
\[
\int_0^1
\rho_N\big(\theta z+(1-\theta)w\big)\,\rd\theta
=
\frac1{|z-w|}
\int_{\min\{z,w\}}^{\max\{z,w\}}\rho_N(s)\,\ds.
\]
Moreover, since $0\le\rho_N\le N$,
\begin{align*}
&\rho_N(z)\,\rho_N(w)\,|z-w|^{\alpha-2}
\int_{\min\{z,w\}}^{\max\{z,w\}}\rho_N(s)\,\ds
\\
&\qquad\le
N^3\chi_{B_R}(z)\,\chi_{B_R}(w)\,|z-w|^{\alpha-1}.
\end{align*}
The function on the right is integrable on $\R^2$, since $\alpha>0$. Thus the resulting
integrand is locally integrable across the diagonal, where it may be defined arbitrarily.
Hence
\[
\Da(\rho_N)
=
\frac12
\int_{\R}\int_{\R}
\rho_N(z) \, \rho_N(w) \, |z-w|^{\alpha-2}
\int_{\min\{z,w\}}^{\max\{z,w\}}\rho_N(s)\,\ds
\,\dz\,\dw.
\]

Since $\rho_N=0$ almost everywhere outside $B_R=(-R,R)$, the preceding double integral may be restricted to $B_R\times B_R$. For $z,w\in B_R$, we have
$|z-w|\le2R$ and, since $\alpha-2<0$,
\[
|z-w|^{\alpha-2}\ge(2R)^{\alpha-2}.
\]
Therefore
\[
\Da(\rho_N)
\ge
\frac12(2R)^{\alpha-2}
\int_{\R}\int_{\R}
\rho_N(z) \, \rho_N(w)
\int_{\min\{z,w\}}^{\max\{z,w\}}\rho_N(s)\,\ds
\,\dz\,\dw.
\]

Now, set
\[
M_N:=\|\rho_N\|_{L^1(\R)}
\]
and define
\[
F_N(s):=\int_{-\infty}^s\rho_N(t)\,\dt.
\]
Since $\rho_N\in L^1(\R)$, the function $F_N$ is absolutely continuous and
$F_N'(s)=\rho_N(s)$ 
for almost every $s\in\R$. By symmetry and Tonelli's theorem,
\begin{align*}
&\int_{\R} \int_{\R}
\rho_N(z) \, \rho_N(w)
\int_{\min\{z,w\}}^{\max\{z,w\}}\rho_N(s)\,\ds
\,\dz\,\dw
\\
&\qquad=
2\iint_{\{w<z\}}
\rho_N(z) \, \rho_N(w)
\int_w^z\rho_N(s)\,\ds
\,\dz\,\dw
\\
&\qquad=
2\int_{\R}
\rho_N(s)
\left(
\int_{-\infty}^s\rho_N(w)\,\dw
\right)
\left(
\int_s^\infty\rho_N(z)\,\dz
\right)
\,\ds
\\
&\qquad=
2\int_{\R}
\rho_N(s)F_N(s)(M_N-F_N(s))\,\ds.
\end{align*}
Using the identity
\[
\frac{\rd}{\ds}
\left(
\frac{M_N}{2}F_N(s)^2-\frac13F_N(s)^3
\right)
=
\rho_N(s)F_N(s)\big(M_N-F_N(s)\big)
\]
for almost every $s\in\R$,  and noting that
\[
F_N(s)\to0
\quad\text{as }s\to-\infty,
\qquad
F_N(s)\to M_N
\quad\text{as }s\to\infty,
\]
we obtain
\begin{align*}
2\int_{\R}
\rho_N(s)F_N(s)\big(M_N-F_N(s)\big)\,\ds
=
2\left(
\frac{M_N^3}{2}-\frac{M_N^3}{3}
\right)=
\frac{M_N^3}{3}.
\end{align*}
It follows that
\[
\Da(\rho_N)
\ge
\frac16(2R)^{\alpha-2}M_N^3.
\]

Combining this estimate with \eqref{eq:one-dimensional-D-monotonicity}, we get
\[
\Da(\rho)
\ge
\frac16(2R)^{\alpha-2}
\|\rho_N\|_{L^1(\R)}^3
\qquad
\text{for every }N\in\N.
\]
Finally, since $\rho_N\uparrow\rho$, monotone convergence gives $\|\rho_N\|_{L^1(\R)}
\uparrow
\|\rho\|_{L^1(\R)}$ in $[0,\infty]$. Since $\Da(\rho)<\infty$, the preceding inequality shows that this
limit is finite. Letting $N\to\infty$, we conclude that
\[
\Da(\rho)
\ge
\frac16(2R)^{\alpha-2}
\|\rho\|_{L^1(\R)}^3.
\]
\end{proof}

\section{Box estimates} \label{sec:box-estimates}

In the forthcoming sections we establish lower bounds for $\Ea(\rho)$
and $\Da(\rho)$ valid for specified classes of densities. The bounds
will follow from estimates for characteristic functions of boxes. By a \emph{box} in~$\R^d$ we mean a set obtained from 
\[
\prod_{j=1}^d[-a_j,a_j],
\qquad a_j>0,
\]
by a rigid motion. We call the numbers $a_1,\dots,a_d$ the half-side lengths of the box.

We record for future reference the exponents appearing in the box estimates of this section, which are
\begin{equation}\label{eq:r_alpha}
r_\alpha
:=
\frac{2}{p_\alpha}
=
\begin{dcases}
\alpha+1,
\qquad &0<\alpha\le1,\\[1mm]
2,
\qquad &1\le\alpha<d,
\end{dcases}
\end{equation}
and
\begin{equation}\label{eq:s_alpha}
s_\alpha
:=
\frac{\beta}{q_\alpha}
=
\begin{dcases}
\alpha+1,
\qquad &0<\alpha\le1,\\[3mm]
\frac{2d+\alpha-1}{d},\qquad&1\le\alpha\le3,\\[3mm]
\frac{2d+2}{d},\qquad &3\le\alpha<d.
\end{dcases}
\end{equation}

\subsection{Rigid motion invariance}
We first note that both functionals are invariant under rigid motions in $\R^d$, that is, maps $W:\R^d\to\R^d$  expressed as
\[
W(x)=Ox+x_0,
\]
where $x_0\in\R^d$ and $O$ is an orthogonal matrix, $O^T O = O O^T = I$. Given a density $\rho$, we set 
\begin{equation} \label{eq:rho_rigid_motion}
\rho_W(x)
:=
\rho\big(W^{-1}(x)\big)
=
\rho\big(O^\top(x-x_0)\big),
\qquad
x\in\R^d.
\end{equation}
The stress tensor $S_\alpha(\rho_W)$ at a point $x \in \R^d$ is given by
\begin{equation}\label{eq:S-rigid-motion}
\begin{split}
S_\alpha(\rho_W) 
& =
\tfrac{1}{2} \int_0^1 \int_{\R^d} \rho \big( W^{-1} ( x + (\theta-1)y) \big)  \rho \big(W^{-1} (x+\theta y)\big) |y|^{\alpha-d-2} \, y\otimes y  \,
\dy \, \rd \theta 
\\
& =\tfrac{1}{2} \int_0^1 \int_{\R^d} \rho\big(W^{-1}(x)+(\theta-1)O^\top y \big) \, \rho\big( W^{-1}(x)+\theta O^\top y\big) \, |y|^{\alpha-d-2} \, y\otimes y \, \dy \, \rd \theta 
\\
& = O\,\Sa(\rho)\big(W^{-1}(x)\big)\,O^\top,
\end{split}
\end{equation}
where in the last step we used the change of variables $z=O^\top y$ together with the orthogonality of $O$.

\begin{lemma}\label{lem:rigid-motion-invariance}
Let $\rho$ be a density and $W:\R^d\to\R^d$ be a rigid motion. Then
\begin{equation}\label{eq:ED-rigid-motion}
\Ea(\rho_W)=\Ea(\rho), \qquad \Da(\rho_W)=\Da(\rho),
\end{equation}
where $\rho_W$ is as in~\eqref{eq:rho_rigid_motion}.
\end{lemma}

\begin{proof}
For the energy, make the change of variables
\[
x=W(\xi),
\qquad
y=W(\zeta).
\]
Since $W$ preserves the Lebesgue measure and
\[
|W(\xi)-W(\zeta)|
=
|O(\xi-\zeta)|
=
|\xi-\zeta|,
\]
we obtain
\begin{align*}
\Ea(\rho_W)
&=
\tfrac1{d-\alpha}
\IR \IR
\rho\big(W^{-1}(x)\big) \, \rho\big(W^{-1}(y)\big) \,
|x-y|^{\alpha-d}
\,\dx\,\dy
\\
&=
\tfrac1{d-\alpha}
\IR \IR
\rho(\xi) \, \rho(\zeta) \,
|\xi-\zeta|^{\alpha-d}
\,\rd\xi\,\rd\zeta
\\
&=
\Ea(\rho).
\end{align*}

Next, taking the determinant on both sides of~\eqref{eq:S-rigid-motion} gives
\[
\det \Sa(\rho_W)(x)
=
\det \Sa(\rho)\big(W^{-1}(x)\big).
\]
Hence
\begin{align*}
\Da(\rho_W)
&=
\IR
\big(\rho_W(x)\det \Sa(\rho_W)(x)\big)^{1/d}\,\dx
\\
&=
\IR
\Big(\rho\big(W^{-1}(x)\big)\det \Sa(\rho)\big(W^{-1}(x)\big)\Big)^{1/d}\,\dx.
\end{align*}
The change of variables $w=W^{-1}(x)$ then yields the desired result.
\end{proof}

For every measurable set $E\subseteq\R^d$, one has
\[
(\chi_E)_W=\chi_{W(E)}.
\]
The preceding lemma therefore allows us to assume in all the box estimates below that the box
is centered at the origin and aligned with the coordinate axes. After relabeling the axes, we may
also assume
\[
0<a_d\le\cdots\le a_1.
\]

\subsection{Energy box estimate}

\begin{lemma}\label{lem:energy-box-estimate}
Let $d\in\N$, $0<\alpha<d$, and $R>0$. There exists $c=c(\alpha,d,R)>0$ such that
\[
\mathcal E_\alpha(\chi_Q)
\ge
c \, |Q|^{r_\alpha}
\]
for every box $Q$ whose half-side lengths are at most $R$, where $r_\alpha$ is defined in~\eqref{eq:r_alpha}.
\end{lemma}

\begin{proof}
By rigid-motion invariance, and after relabeling the coordinate axes, we may assume that
\[
Q=\prod_{j=1}^d[-a_j,a_j],
\qquad
0<a_d\le\cdots\le a_1\le R.
\]

Set
\[
\widetilde Q:=\prod_{j=1}^d[-b_j,b_j],
\qquad
b_j:=\frac{a_j}{2}.
\]
If $x\in\widetilde Q$ and $z\in B_{b_d}$, then, for every
$j=1,\ldots,d$,
\[
|x_j+z_j|
\leq
|x_j|+|z_j|
\leq
b_j+b_d
\leq
2b_j
=
a_j.
\]
Hence $x+z\in Q$. Equivalently,
\[
x+B_{b_d}\subseteq Q,
\qquad \mbox{ for $x\in\widetilde Q$}.
\]
Therefore
\begin{align*}
\mathcal E_\alpha(\chi_Q)
&=
\tfrac1{d-\alpha}
\int_Q\int_Q
|x-y|^{\alpha-d}\,\dy\,\dx
\\
&\ge
c_{\alpha,d}
\int_{\widetilde Q}
\int_{x+B_{b_d}}
|x-y|^{\alpha-d}\,\dy\,\dx
\\
&=
c_{\alpha,d}
\int_{\widetilde Q}
\int_{B_{b_d}}
|z|^{\alpha-d}\,\dz\,\dx
\\
&=
c_{\alpha,d} \, |\widetilde Q| \, b_d^\alpha.
\end{align*}
Since $b_j\leq R$ for every $j=1,\ldots,d$, we have
\[
|\widetilde Q|
=
2^d\prod_{j=1}^d b_j
\leq
2^dR^{d-1}b_d,
\]
and so
\[
b_d\geq c_{d,R}|\widetilde Q|.
\]
It follows that
\[
\mathcal E_\alpha(\chi_Q)
\geq
c_{\alpha,d,R}|\widetilde Q|^{\alpha+1}.
\]
Since $|\widetilde Q|=2^{-d}|Q|$, we conclude that
\[
\mathcal E_\alpha(\chi_Q)
\geq c_{\alpha,d,R}|Q|^{\alpha+1}.
\]

On the other hand, since the half-side lengths of $Q$ are at most $R$, we have
\[
|x-y|\le c_d \, R
\qquad
\text{for all }x,y\in Q.
\]
Because $\alpha-d<0$,
\[
|x-y|^{\alpha-d}\ge c_{\alpha,d,R}
\qquad
\text{for all }x,y\in Q.
\]
Therefore
\[
\mathcal E_\alpha(\chi_Q)
\ge
c_{\alpha,d,R}|Q|^2.
\]

Since every half-side length of $Q$ is at most $R$, we have
\[
|Q|\leq(2R)^d.
\]
Thus, if $p\leq q$, then
\[
|Q|^p
=
|Q|^q|Q|^{p-q}
\geq
(2R)^{d(p-q)}|Q|^q.
\]
Hence, after changing the constant, the estimate with the smaller exponent is the stronger one. Since
\[
r_\alpha=\min\{\alpha+1,2\},
\]
the first estimate proves the claim when $0<\alpha\leq1$, while the second proves it when $1\leq\alpha<d$. At $\alpha=1$, the two exponents coincide.
\end{proof}

\subsection{Determinant box estimates}

We consider the determinant functional $\Da(\rho)$ defined in \eqref{eq:D_functional} and first prove a preparatory lemma for estimating it on characteristic functions of boxes.

\begin{lemma}\label{lem:box-reduction}
Let $d\in \N$ and $0<\alpha<d$, let
$$
\begin{aligned}
Q &=\prod_{j=1}^d[-a_j,a_j],
\qquad
0<a_d\le\cdots\le a_1,
\\
\widetilde Q &:=\prod_{j=1}^d[-b_j,b_j],
\qquad
b_j:=\frac{a_j}{2},
\end{aligned}
$$
and define
\[
A_Q:=
\frac12\int_{\widetilde Q} \, |y|^{\alpha-d-2}\,y\otimes y\,\dy.
\]
Then
\begin{equation}\label{eq:box-reduction}
\Da(\chi_Q)
\ge
|\widetilde Q|\,(\det A_Q)^{1/d}.
\end{equation}
\end{lemma}

\begin{proof}
If $x,y\in\widetilde Q$, then, for every $\theta\in[0,1]$,
\[
x+(\theta-1)y, \, x+\theta y \in Q.
\]
Therefore, for every $x\in\widetilde Q$,
\[
\Sa(\chi_Q)(x)\succeq A_Q,
\]
and hence
\[
\det \Sa(\chi_Q)(x)\ge \det A_Q.
\]
Since $\widetilde Q\subseteq Q$, we have $\chi_Q=1$ on $\widetilde Q$. Thus,
\[
\Da(\chi_Q)
\ge
\int_{\widetilde Q}
\big(\!\det \Sa(\chi_Q)(x)\big)^{1/d}\,\dx
\ge
|\widetilde Q|\,(\det A_Q)^{1/d}.
\]
\end{proof}

We next establish three different lower bounds for $\Da(\chi_Q)$ along the class of boxes. The optimal bound will follow by comparing these estimates.

\begin{lemma}\label{lem:first-box-estimate}
Let $d\in \N$, $0<\alpha<d$ and $R>0$. There exists $c=c(\alpha,d,R)>0$ such that
\begin{equation} \label{eq:first-box-estimate}
\Da(\chi_Q)\ge c\,|Q|^{\alpha+1}
\end{equation}
for every box $Q$ whose half-side lengths are at most $R$.
\end{lemma}

\begin{proof}
By rigid-motion invariance, and after relabeling the coordinate axes, we may assume that
\[
Q=\prod_{j=1}^d[-a_j,a_j],
\qquad
0<a_d\le\cdots\le a_1\le R.
\]
Let $\widetilde Q$ and $A_Q$ be as in Lemma~\ref{lem:box-reduction}. Since
$B_{b_d}\subseteq\widetilde Q$, we have
\[
A_Q
\succeq
\frac12\int_{B_{b_d}}|y|^{\alpha-d-2}\,y\otimes y\,\dy.
\]
By symmetry, the off-diagonal entries of the integral vanish and all diagonal entries are equal.
Moreover, for each $i=1,\dots,d$,
\[
\int_{B_{b_d}}|y|^{\alpha-d-2} \, y_i^2\,\dy
=
\frac1d\int_{B_{b_d}}|y|^{\alpha-d}\,\dy
=
c_{\alpha,d} \, b_d^\alpha.
\]
Hence
\[
\int_{B_{b_d}}|y|^{\alpha-d-2}\,y\otimes y\,\dy
=
c_{\alpha,d} \, b_d^\alpha \, I_d.
\]
It follows that
\[
A_Q\succeq c_{\alpha,d} \,  b_d^\alpha \, I_d,
\]
and therefore 
\[
(\det A_Q)^{1/d}\ge c_{\alpha,d} \, b_d^\alpha.
\]
Applying Lemma~\ref{lem:box-reduction}, we get
\[
\Da(\chi_Q)
\ge
c_{\alpha,d} \, |\widetilde Q|\,b_d^\alpha.
\]
Since $b_j\le R$ for every $j=1,\ldots,d$,
\[
|\widetilde Q|
=
2^d\prod_{j=1}^d b_j
\le
2^dR^{d-1}b_d.
\]
Thus
\[
b_d\ge c_{d,R}|\widetilde Q|,
\]
and hence
\[
\Da(\chi_Q)
\ge
c_{\alpha,d,R}|\widetilde Q|^{\alpha+1}.
\]
Since $|\widetilde Q|=2^{-d}|Q|$, the claim follows.
\end{proof}

\begin{lemma}\label{lem:second-box-estimate}
Let $d\in \N$, $0<\alpha<d$ and $R>0$. There exists $c=c(\alpha,d,R)>0$ such that
\begin{equation} \label{eq:second-box-estimate}
\Da(\chi_Q)\ge c\,|Q|^{\frac{2d+\alpha-1}{d}}
\end{equation}
for every box $Q$ whose half-side lengths are at most $R$.
\end{lemma}

\begin{proof}
By rigid-motion invariance, and after relabeling the coordinate axes, we may assume that
\[
Q=\prod_{j=1}^d[-a_j,a_j],
\qquad
0<a_d\le\cdots\le a_1\le R.
\]
Let $A_Q$ and
$\widetilde Q$ be as in~Lemma~\ref{lem:box-reduction}. For each $j=1,\dots,d$, we write
\begin{equation} \label{eq:eigenvalue_A_Q}
\lambda_j:= (A_Q)_{jj}
=
\frac12\int_{\widetilde Q}|y|^{\alpha-d-2} \, y_j^2\,\dy,
\end{equation}
and set
\[
E_j:=
\prod_{i=1}^{j}[-b_j,b_j]
\times
\prod_{i=j+1}^{d}[-b_i,b_i], \qquad b_j = \frac{a_j}{2}.
\]
Since $b_1\ge\cdots\ge b_d$, we have $E_j\subseteq\widetilde Q$. Moreover, for $y\in E_j$,
$|y|\le\sqrt d\,b_j$ and so
\[
|y|^{\alpha-d-2}\ge c_{\alpha,d} \, b_j^{\alpha-d-2},
\qquad y\in E_j.
\]
Hence
\[
\lambda_j
\ge
c_{\alpha,d} \, b_j^{\alpha-d-2}
\int_{E_j}y_j^2\,\dy.
\]
Since
\begin{align*}
\int_{E_j}y_j^2\,\dy &= \bigg(\prod_{i=1}^{j-1} \int_{-b_j}^{b_j} 1 \, \dy_i \bigg)\bigg( \int_{-b_j}^{b_j} y_j^2 \, \dy_j\bigg)\bigg(\prod_{i=j+1}^{d} \int_{-b_i}^{b_i} 1 \, \dy_i \bigg) \\
&=c_d\,b_j^{j+2}\prod_{i=j+1}^d b_i,
\end{align*}
we obtain
\begin{equation}\label{eq:lambda-box-estimate}
\lambda_j
\ge
c_{\alpha,d} \, b_j^{\alpha-d+j}\prod_{i=j+1}^d b_i.
\end{equation}
Multiplying \eqref{eq:lambda-box-estimate} over $j=1,\dots,d$ and using that $A_Q$ is diagonal,
we get
\[
\det A_Q
=
\prod_{j=1}^d\lambda_j
\ge
c_{\alpha,d}
\prod_{j=1}^d
\Big(
b_j^{\alpha-d+j}\prod_{i=j+1}^d b_i
\Big).
\]
For each fixed $k$, the factor $b_k$ appears first in the term $j=k$, with exponent $\alpha-d+k$, and then once in each product with $j<k$. Hence it appears $k-1$ more times, and so the total exponent of $b_k$ is
\[
\alpha-d+k+(k-1)=\alpha-d+2k-1.
\]
Therefore
\[
\det A_Q
\ge
c_{\alpha,d}
\prod_{k=1}^d b_k^{\alpha-d+2k-1},
\]
and hence
\[
(\det A_Q)^{1/d}
\ge
c_{\alpha,d}
\prod_{k=1}^d b_k^{\frac{\alpha-d+2k-1}{d}}.
\]
Using Lemma~\ref{lem:box-reduction} and $|\widetilde Q|=2^d\prod_{k=1}^d b_k$, 
we get
\[
\Da(\chi_Q)
\ge
c_{\alpha,d}
\prod_{k=1}^d b_k^{1+\frac{\alpha-d+2k-1}{d}}
=
c_{\alpha,d}
\prod_{k=1}^d b_k^{\frac{\alpha+2k-1}{d}}.
\]
Since
\[
\frac{\alpha+2k-1}{d}
\le
\frac{2d+\alpha-1}{d}, \qquad \text{and} \qquad  b_k \leq R, \quad k=1,\ldots,d,
\]
we have
\[
b_k^{\frac{\alpha+2k-1}{d}}
\ge
\Big(\min_{\ell=1,\ldots,d} R^{\frac2d(\ell-d)} \Big) \, 
b_k^{\frac{2d+\alpha-1}{d}}.
\]
Thus
\[
\Da(\chi_Q)
\ge
c_{\alpha,d,R}
\Big(\prod_{k=1}^d b_k\Big)^{\frac{2d+\alpha-1}{d}}.
\]
The claim follows by noting that $|Q|=4^d\prod_{k=1}^d b_k$.
\end{proof}

\begin{lemma}\label{lem:third-box-estimate}
Let $d\in \N$, $0<\alpha<d$ and $R>0$. There exists $c=c(\alpha,d,R)>0$ such that
\begin{equation}\label{eq:third-box-estimate}
\Da(\chi_Q)
\ge
c\,|Q|^{\frac{2d+2}{d}}
\end{equation}
for every box $Q$ whose half-side lengths are at most $R$.
\end{lemma}

\begin{proof}
By rigid-motion invariance, and after relabeling the coordinate axes, we may assume that
\[
Q=\prod_{j=1}^d[-a_j,a_j],
\qquad
0<a_d\le\cdots\le a_1\le R.
\]
Let $A_Q$ and
$\widetilde Q$ be as in~Lemma~\ref{lem:box-reduction}, and let $\lambda_j$ be as in~\eqref{eq:eigenvalue_A_Q}.

Since $b_j\le R/2$ for every $j$, we have $|y|\le c_d \, R$ for $y\in\widetilde Q$ and thus
\[
|y|^{\alpha-d-2}\ge c_{\alpha,d,R}
\qquad\text{for all }y\in\widetilde Q.
\]
Therefore, \eqref{eq:eigenvalue_A_Q} gives
\[
\lambda_j
\ge
c_{\alpha,d,R}
\int_{\widetilde Q}y_j^2\,\dy = c_{\alpha,d,R} \,
b_j^2\prod_{i=1}^d b_i.
\]
Multiplying over $j=1,\dots,d$, we obtain
\[
\det A_Q
=
\prod_{j=1}^d\lambda_j
\ge
c_{\alpha,d,R}
\prod_{j=1}^d b_j^2 \, \Big(\prod_{i=1}^d b_i\Big)^d
=
c_{\alpha,d,R}
\Big(\prod_{i=1}^d b_i\Big)^{d+2}.
\]
Thus
\[
(\det A_Q)^{1/d}
\ge
c_{\alpha,d,R}
\Big(\prod_{i=1}^d b_i\Big)^{1+\frac2d}.
\]
Using Lemma~\ref{lem:box-reduction}, we get
\[
\Da(\chi_Q)
\ge
c_{\alpha,d,R}
\Big(\prod_{i=1}^d b_i\Big)^{2+\frac2d}.
\]
Finally, $|Q|=4^d\prod_{i=1}^d b_i$, and therefore
\[
\Da(\chi_Q)
\ge
c_{\alpha,d,R} \, |Q|^{2+\frac2d}.
\]
This proves \eqref{eq:third-box-estimate}.
\end{proof}

The three preceding estimates yield the following result.

\begin{lemma}\label{lem:combined-determinant-box-estimate}
Let $d\in\N$, $0<\alpha<d$, and $R>0$. There exists
$c=c(\alpha,d,R)>0$ such that
\[
\Da(\chi_Q)\geq c\,|Q|^{s_\alpha}
\]
for every box $Q$ whose half-side lengths are at most $R$, where
$s_\alpha$ is defined in~\eqref{eq:s_alpha}.
\end{lemma}

\begin{proof}
The preceding three lemmas give the exponents
\[
p_1:=\alpha+1,
\qquad
p_2:=\frac{2d+\alpha-1}{d},
\qquad
p_3:=\frac{2d+2}{d}.
\]
Since $|Q|\leq(2R)^d$, after adjusting the constant, the estimate
with the smallest exponent is the strongest. Moreover,
\[
p_2-p_1
=
\frac{(d-1)(1-\alpha)}{d},
\qquad
p_3-p_2
=
\frac{3-\alpha}{d}.
\]
It follows that the smallest exponent is $p_1$ when
$0<\alpha\leq1$, $p_2$ when $1\leq\alpha\leq3$, and $p_3$ when $3\leq\alpha<d$. Thus the smallest exponent is $s_\alpha$,
as defined in~\eqref{eq:s_alpha}.
\end{proof}

\section{Lorentz and Lebesgue estimates} \label{sec:lorentz-quasiconcave}

In this section we prove Theorem~\ref{thm:ED-log-concave}. First, using John's ellipsoid theorem,  the box estimates of Section~\ref{sec:box-estimates} are extended to characteristic functions of
convex sets. The latter estimates are applied to the superlevel sets of quasiconcave densities and, after integration over the level parameter, provide endpoint Lorentz bounds. A weak-to-strong lemma for log-concave densities then yields the corresponding strong
Lebesgue estimates.

\subsection{Estimates for convex sets} \label{sec:convex-lower-bounds}
We start with the estimates on characteristic functions $\chi_K$ of convex, bounded sets.

\begin{proposition}\label{prop:ED_convex}
Let $d\in\N$, $0<\alpha<d$, and $R>0$. There exists
$c=c(\alpha,d,R)>0$ such that
\[
\Ea(\chi_K)\geq c\,|K|^{r_\alpha}
\qquad \text{and} \qquad
\Da(\chi_K)\geq c\,|K|^{s_\alpha}
\]
for every convex set $K\subseteq B_R$, where $r_\alpha$ and
$s_\alpha$ are defined in~\eqref{eq:r_alpha} and~\eqref{eq:s_alpha},
respectively.
\end{proposition}

\begin{proof}
If $|K|=0$, then $\chi_K=0$ a.e., and both estimates
are trivial. We may therefore assume that $|K|>0$. Since $K$ is convex and has positive measure, it has nonempty interior. Moreover, its boundary has Lebesgue measure zero. Hence $|K|=|\overline K|$ and
$\chi_K=\chi_{\overline K} \; \text{a.e. in }\R^d$. 
Replacing $K$ by its closure, we may assume that $K$ is a convex
body contained in $\overline{B_R}$.

By John's ellipsoid theorem~\cite[Theorem~10.12.2]{schneider2014convex}, there exists
an ellipsoid $E$ with center~$c_E$ such that
\[
E\subseteq K\subseteq c_E+d(E-c_E),
\]
and so,
\[
|E|\geq d^{-d}|K|.
\]

Now, we write $E$ in principal-axis coordinates as
\[
E
=
c_E+
\bigg\{
\sum_{j=1}^d t_jv_j:\,
\sum_{j=1}^d\frac{t_j^2}{\ell_j^2}\leq1
\bigg\},
\]
where $v_1,\ldots,v_d$ form an orthonormal basis of principal directions of $E$, and $\ell_1,\ldots,\ell_d>0$ are the corresponding semiaxis lengths. Let $Q$ be the inscribed box
\[
Q
:=
c_E+
\bigg\{
\sum_{j=1}^d t_jv_j:\,
|t_j|\leq\frac{\ell_j}{\sqrt d},
\quad j=1,\ldots,d
\bigg\}.
\]
Then $Q\subseteq E$, since
\[
\sum_{j=1}^d\frac{t_j^2}{\ell_j^2}
\leq
\sum_{j=1}^d\frac1d
=
1.
\]
Moreover,
\[
|Q|
=
\bigg(\frac2{\sqrt d}\bigg)^d
\prod_{j=1}^d\ell_j
=
c_d|E|
\geq
c_d|K|.
\]

Since $Q\subseteq E \subseteq K \subseteq \overline{B_R}$, every half-side length of $Q$ is at most $R$. Applying Lemma~\ref{lem:energy-box-estimate} to $Q$ gives
\[
\Ea(\chi_K)
\geq
\Ea(\chi_Q)
\geq
c_{\alpha,d,R}|Q|^{r_\alpha}
\geq
c_{\alpha,d,R}|K|^{r_\alpha}.
\]
Likewise, $\chi_Q\leq\chi_K$ implies
\[
\Sa(\chi_Q)(x)\preceq\Sa(\chi_K)(x),
\qquad x\in\R^d.
\]
By the monotonicity of  determinant functions along the L\"{o}wner order,
\[
\det\Sa(\chi_Q)(x)
\leq
\det\Sa(\chi_K)(x),
\]
and therefore
\begin{align*}
\Da(\chi_K)
&\geq
\int_Q
\big(\!\det\Sa(\chi_K)(x)\big)^{1/d}\,\dx\\
&\geq\!
\int_Q
\big(\!\det\Sa(\chi_Q)(x)\big)^{1/d}\,\dx\\
&=
\Da(\chi_Q).
\end{align*}
Using Lemma~\ref{lem:combined-determinant-box-estimate}, we obtain
\[
\Da(\chi_K)
\geq
c_{\alpha,d,R}|Q|^{s_\alpha}
\geq
c_{\alpha,d,R}|K|^{s_\alpha}.
\]
\end{proof}

\begin{remark*}
For $d\ge2$ and $0<\alpha\le1$, the energy estimate in Proposition~\ref{prop:ED_convex} can alternatively be deduced from an interpolation inequality of Kn\"upfer and Muratov~\cite[Lemma~7.1]{knupfer2014isoperimetric};
see also~\cite{knupfer2013isoperimetric}.
For a convex body $K$, their inequality gives, in our notation,
\[
\Ea(\chi_K)
\geq c_{\alpha,d}\frac{|K|^{1+\alpha}}{P(K)^\alpha},
\]
where $P(K)$ denotes its perimeter. Since perimeter is monotone under inclusion for convex bodies (see~\cite{stefani2018monotonicity}) and
$K\subseteq B_R$, we have $P(K)\leq P(B_R)$,
yielding the exponent $1+\alpha=r_\alpha$. For $1\le\alpha<d$, the confined $L^1$ bound in Proposition~\ref{prop:E_L1_confined} gives the exponent $r_\alpha=2$, improving this volume bound when $\alpha>1$.
\end{remark*}

\subsection{Lorentz estimates for quasiconcave densities}

A density $\rho$ is called \emph{quasiconcave} if
\begin{equation}\label{eq:quasiconcavity}
\rho\big((1-\theta)x+\theta y\big)
\geq
\min\{\rho(x),\rho(y)\}
\end{equation}
for every $x,y\in\R^d$ and every $\theta\in[0,1]$.
Equivalently, $\rho$ is quasiconcave if each superlevel set
$\{\rho>t\}$, with $t>0$, is convex; see~\cite[Chapter~8]{sundaram1996first}.

In the sequel we use the following distribution-function formulation of the
Lorentz quasi-norms. For a density $\rho$ and indices
$0<p<\infty$ and $0<s\leq\infty$, we write
\begin{equation}\label{eq:Lorentz_norm}
\|\rho\|_{L^{p,s}(\R^d)}
=
\begin{dcases}
\displaystyle
\bigg(
s\int_0^\infty
t^{s-1}\,\big|\{\rho>t\}\big|^{s/p}\,\dt
\bigg)^{1/s},
&\qquad 0<s<\infty,\\[2mm]
\displaystyle
\sup_{t>0}
t\,\big|\{\rho>t\}\big|^{1/p},
&\qquad s=\infty.
\end{dcases}
\end{equation}
We adopt the convention $L^{\infty,\infty}(\R^d)=L^\infty(\R^d)$.
The Lorentz quasi-norm for $s=\infty$ is also known as the weak $L^p$ quasi-norm. For $0<p<\infty$, we recall that the $L^p$ quasi-norm agrees with the $L^{p,p}$ Lorentz quasi-norm,
\[
\|\rho\|_{L^p(\R^d)}
=
\bigg(
p\int_0^\infty
t^{p-1}\,\big|\{\rho>t\}\big|\,\dt
\bigg)^{1/p},
\]
and that the continuous embeddings $L^{p,q}(\R^d)\subseteq L^{p,r}(\R^d)$ hold for $0<q\le r\le\infty$.

When $\rho$ is supported in $B_R$, we write
$\|\rho\|_{L^{p,s}(B_R)}$ for the expression
in~\eqref{eq:Lorentz_norm} with $\{\rho>t\}$ replaced by
$\{\rho>t\}\cap B_R$. The formulation
in~\eqref{eq:Lorentz_norm} is equivalent to the usual definition
in terms of decreasing rearrangements;
see~\cite[Chapter~1]{grafakos2014classical}.

\begin{lemma}\label{lem:layer-cake}
Let $d\in\N$, $0<\alpha<d$, and $R>0$. Let $\rho$ be a density
supported in $B_R$, and, for each $t>0$, set
\[
E_t:=\{\rho>t\}\cap B_R.
\]
Then
\begin{equation}\label{eq:energy-layer-cake}
\Ea(\rho)
\geq
2\int_0^\infty
t\,\Ea(\chi_{E_t})\,\dt
\end{equation}
and
\begin{equation}\label{eq:determinant-layer-cake}
\Da(\rho)
\geq
2\int_0^\infty
t^{\beta-1}\Da(\chi_{E_t})\,\dt.
\end{equation}
\end{lemma}

\begin{proof}
For $N\in\N$, set $\rho_N:=\min\{\rho,N\}$. One easily checks that
$$
\begin{aligned}
\Ea(\rho)  &\geq \Ea(\rho_N),
\\
\Sa(\rho)  \succeq \Sa(\rho_N) \quad & \qquad \det \Sa(\rho) \ge \det \Sa(\rho_N)
\\
\Da(\rho) &\geq\Da(\rho_N).
\end{aligned}
$$

For nonnegative numbers $a,b$, we have
\begin{equation} \label{eq:aux_inequality_quasiconcave}
\min\{a,N\}\min\{b,N\}
\geq
\min\{a,b,N\}^2
=
\int_0^N
2t\,\chi_{\{a>t\}} \, \chi_{\{b>t\}}\,\dt.
\end{equation}
Applying this inequality with $a=\rho(x)$, $b=\rho(y)$ and interchanging the order of integration, we obtain
\[
\Ea(\rho)
\geq
\Ea(\rho_N)
\geq
2\int_0^N
t\,\Ea(\chi_{E_t})\,\dt.
\]
As the integrand is nonnegative, taking the supremum over~$N\in\N$ gives~\eqref{eq:energy-layer-cake}.

To prove the determinant estimate, apply the inequality~\eqref{eq:aux_inequality_quasiconcave} with
$a=\rho\big(x+(\theta-1)y\big)$, $b=\rho\big(x+\theta y\big)$ to obtain
\[
\Sa(\rho_N)(x)
\succeq
\int_0^N
2t\,\Sa(\chi_{E_t})(x)\,\dt.
\]
Let $\Phi(A):=(\det A)^{1/d}$ be defined on the cone of positive semidefinite $d\times d$ matrices. By
Minkowski's determinant inequality~\cite[Theorem~7.8.21]{horn2013matrix}, $\Phi$ is concave. Using Jensen's inequality we deduce
\begin{align*}
\Phi\big(\Sa(\rho_N)(x)\big)
&\geq
\Phi\left(
\int_0^N
2t\,\Sa(\chi_{E_t})(x)\,\dt
\right)
\\
&\geq
\int_0^N
2t\,
\Phi\big(\Sa(\chi_{E_t})(x)\big)\, \dt.
\end{align*}
Therefore, using Tonelli's theorem and then restricting the inner
integral to $E_t$, we get
\begin{align*}
\Da(\rho)
&\geq
\Da(\rho_N)
=
\int_{\R^d}
\rho_N(x)^{1/d}\,
\Phi\big(\Sa(\rho_N)(x)\big)\,\dx
\\
&\geq
2\int_0^N
t
\int_{\R^d}
\rho_N(x)^{1/d}\,
\Phi\big(\Sa(\chi_{E_t})(x)\big)
\,\dx\,\dt
\\
&\geq
2\int_0^N
t
\int_{E_t}
\rho_N(x)^{1/d}\,
\Phi\big(\Sa(\chi_{E_t})(x)\big)
\,\dx\,\dt.
\end{align*}

For $0<t<N$, we have $\rho_N(x)>t$ for $x\in E_t$. It follows
\begin{align*}
\Da(\rho)
&\geq
2\int_0^N
t^{1+\frac1d}
\int_{E_t}
\Phi\big(\Sa(\chi_{E_t})(x)\big)
\,\dx\,\dt
\\
&=
2\int_0^N
t^{\beta-1}\Da(\chi_{E_t})\,\dt
\end{align*}
and letting $N\to\infty$ proves
\eqref{eq:determinant-layer-cake}.
\end{proof}

Combining Proposition~\ref{prop:ED_convex} with Lemma~\ref{lem:layer-cake} gives the following endpoint
Lorentz estimates.

\begin{theorem}\label{thm:ED-quasiconcave}
Let $d\in\N$, $0<\alpha<d$, and $R>0$. There exists a
positive constant $C=C(\alpha,d,R)$ such that
\begin{equation}\label{eq:ED-quasiconcave}
\Ea(\rho)
\geq
C\|\rho\|_{L^{p_\alpha,2}(B_R)}^2
\qquad\text{and}\qquad
\Da(\rho)
\geq
C\|\rho\|_{L^{q_\alpha,\beta}(B_R)}^\beta
\end{equation}
for all quasiconcave densities $\rho$ supported in $B_R$.
\end{theorem}

\begin{proof}
Let $\rho$ be a quasiconcave density supported in $B_R$, and, for each
$t>0$, set
\[
E_t:=\{\rho>t\}\cap B_R.
\]
Since $\{\rho>t\}$ and $B_R$ are convex, $E_t$ is a convex set
contained in $B_R$. Proposition~\ref{prop:ED_convex} therefore gives
\[
\Ea(\chi_{E_t})
\geq
c_{\alpha,d,R}\,|E_t|^{r_\alpha},
\qquad
\Da(\chi_{E_t})
\geq
c_{\alpha,d,R}\,|E_t|^{s_\alpha}.
\]

For the energy functional, using~\eqref{eq:energy-layer-cake}, we obtain
\begin{equation}\label{eq:energy-layers}
\Ea(\rho)
\geq
 c_{\alpha,d,R}
\int_0^\infty
t\,|E_t|^{r_\alpha}\,\dt.
\end{equation}
Since $r_\alpha=2/p_\alpha$, the Lorentz quasi-norm becomes 
\[
\|\rho\|_{L^{p_\alpha,2}(B_R)}^2
=
2\int_0^\infty
t\,|E_t|^{r_\alpha}\,\dt , 
\]
and \eqref{eq:energy-layers} yields the estimate for the energy functional in~\eqref{eq:ED-quasiconcave}.

For the determinant functional, using~\eqref{eq:determinant-layer-cake}, we get
\begin{equation}\label{eq:determinant-layers}
\Da(\rho)
\geq
2c_{\alpha,d,R}
\int_0^\infty
t^{\beta-1}|E_t|^{s_\alpha}\,\dt.
\end{equation}
Since $s_\alpha=\beta/q_\alpha$, 
\[
\|\rho\|_{L^{q_\alpha,\beta}(B_R)}^\beta =
\beta
\int_0^\infty
t^{\beta-1}|E_t|^{s_\alpha}\,\dt,
\]
and \eqref{eq:determinant-layers} yields the estimate for the determinant functional in~\eqref{eq:ED-quasiconcave}.
\end{proof}

\subsection{Weak Lebesgue estimates}

The endpoint estimates proved above yield the corresponding weak endpoint
bounds. Since the densities are supported in the finite-measure set $B_R$,
they also yield lower bounds in Lorentz spaces with every smaller first
index.

\begin{corollary}\label{cor:quasiconcave_weak_lebesgue}
Let $d\in\N$, $0<\alpha<d$, and $R>0$. There exists
$C=C(\alpha,d,R)>0$ such that
\[
\Ea(\rho)
\geq
C \, \|\rho\|_{L^{p_\alpha,\infty}(B_R)}^2 \qquad \text{and} \qquad
\Da(\rho)
\geq
C \, \|\rho\|_{L^{q_\alpha,\infty}(B_R)}^\beta
\]
for every quasiconcave density $\rho$ supported in $B_R$, where $p_\alpha$, $q_\alpha$, and $\beta$ are defined in~\eqref{eq:p_alpha}, \eqref{eq:q_alpha}, and~\eqref{eq:beta}, respectively.

Moreover, the following statements hold.

\begin{enumerate}
\item[(i)]
For every $0<p<p_\alpha$ and every $0<s\leq\infty$, there exists
$C=C(\alpha,d,p,s,R)>0$ such that
\[
\Ea(\rho)
\geq
C \, \|\rho\|_{L^{p,s}(B_R)}^2
\]
for every quasiconcave density $\rho$ supported in $B_R$.

\item[(ii)]
For every $0<q<q_\alpha$ and every $0<s\leq\infty$, there exists
$C=C(\alpha,d,q,s,R)>0$ such that
\[
\Da(\rho)
\geq
C \, \|\rho\|_{L^{q,s}(B_R)}^\beta
\]
for every quasiconcave density $\rho$ supported in $B_R$.
\end{enumerate}
\end{corollary}

\begin{proof}
The continuous Lorentz embeddings
$
L^{p_\alpha,2}(\R^d)
\subseteq
L^{p_\alpha,\infty}(\R^d)
$
and
$
L^{q_\alpha,\beta}(\R^d)
\subseteq
L^{q_\alpha,\infty}(\R^d),
$
together with Theorem~\ref{thm:ED-quasiconcave}, give the two weak endpoint estimates.

Since $|B_R|<\infty$, for every $0<r<a<\infty$ and
$0<s\leq\infty$, there exists $C=C(a,d,r,s,R)>0$ such that
\[
\|\rho\|_{L^{r,s}(B_R)}
\leq
C\|\rho\|_{L^{a,\infty}(B_R)}.
\]
Indeed, for $0 < s < \infty$, this follows directly from
\[
\big|\{\rho>t\}\cap B_R\big|
\leq
\min\big\{
|B_R|, \,t^{-a}
\|\rho\|_{L^{a,\infty}(B_R)}^a 
\big\}
\]
by splitting the distribution-function integral at
$t=\|\rho\|_{L^{a,\infty}(B_R)}|B_R|^{-1/a}$. For $s=\infty$, the same distribution-function bound gives the estimate directly. Applying this embedding with $a=p_\alpha$ and $r=p$ gives part~(i), while
applying it with $a=q_\alpha$ and $r=q$ gives part~(ii).
\end{proof}

In particular, taking $s=p$ in part~{\rm (i)} and $s=q$ in
part~{\rm (ii)}, we obtain
\[
\Ea(\rho)
\geq
C\|\rho\|_{L^p(B_R)}^2,
\qquad
0<p<p_\alpha,
\]
and
\[
\Da(\rho)
\geq
C\|\rho\|_{L^q(B_R)}^\beta,
\qquad
0<q<q_\alpha.
\]
In summary, for densities supported in a fixed ball that are quasiconcave, strong $L^p$ and
$L^q$ lower bounds are valid for every index strictly less than the corresponding
endpoint. The strong endpoint estimates under log-concavity are proved in
the next subsection.

\subsection{Lebesgue estimates for log-concave densities} \label{sec:lebesgue-log-concave}

We first prove that the weak $L^q$ quasi-norm controls the strong
$L^q$ quasi-norm within the log-concave class.

\begin{lemma}\label{lem:log-concave-weak-strong}
Let $d \in \N$ and $0<q<\infty$. There exists $c=c(d,q)>0$ such that
\[
\|\rho\|_{L^q(\R^d)}
\leq
c \,\|\rho\|_{L^{q,\infty}(\R^d)}
\]
for all log-concave densities $\rho$.
\end{lemma}

\begin{proof}
It suffices to prove the estimate for bounded densities. Indeed, for a
general log-concave density $\rho$, the truncations
$\rho_N:=\min\{\rho,N\}$, 
$N\in\N$, 
are log-concave, and the result follows from the bounded case by monotone convergence.

If $\rho=0$ almost everywhere or $\|\rho\|_{L^{q,\infty}(\R^d)}=\infty$, there is nothing to prove. Set
\[
W:=\|\rho\|_{L^{q,\infty}(\R^d)},
\qquad
M:=\|\rho\|_{L^\infty(\R^d)},
\]
and assume $0<M,W<\infty$. Choose $x_0\in\R^d$ such that
$\rho(x_0)>\frac M2$ and for each $t>0$, let $E_t:=\{\rho>t\}$. 
If $0<t<M$ and $x\in E_t$, then log-concavity gives
\[
\rho\left(\frac{x_0+x}{2}\right)
\geq
\rho(x_0)^{1/2}\rho(x)^{1/2}>
\left(\frac{Mt}{2}\right)^{1/2},
\]
and so,
\[
\frac{x_0 + E_t}{2}
\subseteq
E_{\sqrt{Mt/2}}.
\]
Taking Lebesgue measures and using
the definition of $W$, we obtain
\[
2^{-d}|E_t|
\leq
\left|E_{\sqrt{Mt/2}}\right|\leq
W^q\left(\frac{Mt}{2}\right)^{-q/2}.
\]
Thus
\[
|E_t|
\leq
2^{d+q/2}W^qM^{-q/2}t^{-q/2},
\qquad
0<t<M.
\]

Consequently,
\begin{align*}
\|\rho\|_{L^q(\R^d)}^q
&=
q\int_0^M t^{q-1}|E_t|\,\dt\\
&\leq
q\,2^{d+q/2}W^qM^{-q/2}
\int_0^M t^{q/2-1}\,\dt\\
&=
2^{d+1+q/2}W^q,
\end{align*}
which yields the desired result.
\end{proof}

\begin{proof}[Proof of Theorem~\ref{thm:ED-log-concave}]
Let $\rho$ be a log-concave density supported in $B_R$. Since every
log-concave density is quasiconcave,
Corollary~\ref{cor:quasiconcave_weak_lebesgue} gives
\[
\Ea(\rho)
\geq
C_{\alpha,d,R}
\|\rho\|_{L^{p_\alpha,\infty}(B_R)}^2,
\qquad
\Da(\rho)
\geq
C_{\alpha,d,R}
\|\rho\|_{L^{q_\alpha,\infty}(B_R)}^\beta.
\]

On the other hand, since $\rho$ is supported in $B_R$, Lemma~\ref{lem:log-concave-weak-strong},
applied with $q=p_\alpha$ and $q=q_\alpha$, gives
\[
\|\rho\|_{L^{p_\alpha}(B_R)}
\leq
C_{\alpha,d}
\|\rho\|_{L^{p_\alpha,\infty}(B_R)},
\qquad
\|\rho\|_{L^{q_\alpha}(B_R)}
\leq
C_{\alpha,d}
\|\rho\|_{L^{q_\alpha,\infty}(B_R)}.
\]

Combining these estimates yields~\eqref{eq:ED-log-concave} and concludes the proof.
\end{proof}

\section{Sharpness of the critical exponents} \label{sec:sharpness}

In this section we prove Theorem~\ref{thm:sharpness_1}.
Using normalized characteristic functions of thin plates, we show that the integrability exponents $p_\alpha$ and $q_\alpha$ in Theorem~\ref{thm:ED-log-concave} are optimal.

\subsection{Thin-plate estimates}

For $0<\eps<1$, consider the set
\begin{equation}\label{eq:sharpness-thin-plate}
T_\eps
:=
\begin{dcases}
(-\eps,\eps),
\qquad & d=1,\\[1mm]
B_1^{d-1}\times(-\eps,\eps),
\qquad & d\geq2,
\end{dcases}
\end{equation}
and for $0<r<\infty$, define
\begin{equation}\label{eq:sharpness-normalized-thin-plate}
\rho_{\eps,r}
:=
|T_\eps|^{-1/r}\chi_{T_\eps}.
\end{equation}
The set $T_\eps$ is convex, so $\rho_{\eps,r}$ is log-concave. Its support satisfies $\supp\rho_{\eps,r}\subseteq B_{\sqrt2}$, and for $0<s<\infty$ one easily computes
\begin{align*}
\|\rho_{\eps,r}\|_{L^{r,s}(\R^d)}^s
&=
s\int_0^{|T_\eps|^{-1/r}}
t^{s-1}|T_\eps|^{s/r}\,\dt=
1,
\end{align*}
and
\[
\|\rho_{\eps,r}\|_{L^{r,\infty}(\R^d)}
=
\sup_{0<t<|T_\eps|^{-1/r}}
t|T_\eps|^{1/r}
=
1.
\]
In summary, the Lorentz quasi-norm is normalized,
\begin{equation}\label{eq:thin-plate-lorentz}
\|\rho_{\eps,r}\|_{L^{r,s}(\R^d)}=1,
\qquad
\text{for }0<s\leq\infty.
\end{equation}
We proceed to estimate the energy $\Ea$ and determinant $\Da$ along these examples.

\begin{lemma}\label{lem:energy-thin-plate-upper-bound}
Let $d\in\N$ and $0<\alpha<d$. Then
\begin{equation}\label{eq:estim-energy}
\Ea(\chi_{T_\eps})
\leq
\begin{dcases}
c_{\alpha,d}\,\eps^{\alpha+1},
\qquad &0<\alpha<1,\\[1mm]
c_d\,\eps^2\big(1+|\!\log\eps|\big),
\qquad &\alpha=1,\\[1mm]
c_{\alpha,d}\,\eps^2,
\qquad &1<\alpha<d.
\end{dcases}
\end{equation}
\end{lemma}

\begin{proof}
If $d=1$, then $0<\alpha<1$, and the change of variables $x=\eps x'$ and $y=\eps y'$ gives
\[
\Ea(\chi_{T_\eps})
=
\eps^{\alpha+1}\Ea(\chi_{(-1,1)})
=
c_\alpha\,\eps^{\alpha+1}.
\]

Assume now that $d\geq2$, and write
$x=(u,t)$ and $y=(v,s)$, where $u,v\in\R^{d-1}$ and
$t,s\in\R$. By Tonelli's theorem, followed by the change of variables $t=\eps\tau$ and $s=\eps\sigma$, we obtain
\begin{equation} \label{eq:thin-plates-energy}
\begin{split}
\Ea(\chi_{T_\eps})
&=
\frac1{d-\alpha}
\int_{-\eps}^{\eps}
\int_{-\eps}^{\eps}
G_{\alpha,d}(t-s)\,\dt\,\ds\\
&=
\frac{\eps^2}{d-\alpha}
\int_{-1}^1
\int_{-1}^1
G_{\alpha,d}\big(\eps(\tau-\sigma)\big)
\,\rd\tau\,\rd\sigma,
\end{split}
\end{equation}
where
\begin{equation}\label{eq:G_alpha,d}
G_{\alpha,d}(h)
:=
\int_{B_1^{d-1}}
\int_{B_1^{d-1}}
\big(h^2+|u-v|^2\big)^{\frac{\alpha-d}{2}}
\,\du\,\dv.
\end{equation}

We estimate $G_{\alpha,d}(h)$ for $0<|h|\leq2$.
Making the change of variables $\eta=u-v$, using polar
coordinates in $\R^{d-1}$, and splitting the radial integral
at $r=|h|$, we get
\begin{align*}
G_{\alpha,d}(h)
&\leq
c_d
\int_{|\eta|\leq2}
\big(h^2+|\eta|^2\big)^{\frac{\alpha-d}{2}}
\,\rd\eta\\
&\leq
c_d\bigg(
|h|^{\alpha-d}
\int_0^{|h|}r^{d-2}\,\rd r+\int_{|h|}^2r^{\alpha-2}\,\rd r\bigg).
\end{align*}
The first term inside the parentheses equals
$|h|^{\alpha-1}/(d-1)$. Moreover,
\[
\int_{|h|}^2r^{\alpha-2}\,\rd r
=
\begin{dcases}
\frac{|h|^{\alpha-1}-2^{\alpha-1}}{1-\alpha},
&\qquad 0<\alpha<1,\\[3mm]
\displaystyle
\log\frac2{|h|},
&\qquad \alpha=1,\\[3mm]
\displaystyle
\frac{2^{\alpha-1}-|h|^{\alpha-1}}{\alpha-1},
&\qquad 1<\alpha<d.
\end{dcases}
\]
Consequently, for $0<|h|\leq2$,
\begin{equation}\label{eq:sharpness-energy-effective-kernel}
G_{\alpha,d}(h)
\leq
\begin{dcases}
c_{\alpha,d}\,|h|^{\alpha-1},
&\qquad 0<\alpha<1,\\[1mm]
c_d\big(1+\big|\!\log|h|\big|\big),
&\qquad \alpha=1,\\[1mm]
c_{\alpha,d},
&\qquad 1<\alpha<d.
\end{dcases}
\end{equation}

Returning to the energy identity~\eqref{eq:thin-plates-energy}, if $0<\alpha<1$,
then
\begin{align*}
\Ea(\chi_{T_\eps})
&\leq
c_{\alpha,d}\,\eps^{\alpha+1}
\int_{-1}^1
\int_{-1}^1
|\tau-\sigma|^{\alpha-1}
\,\rd\tau\,\rd\sigma\\
&\leq
c_{\alpha,d}\,\eps^{\alpha+1}.
\end{align*}
If $\alpha=1$, then
\begin{align*}
\Ea(\chi_{T_\eps})
&\leq
c_d\,\eps^2
\int_{-1}^1
\int_{-1}^1
\big(
1+|\!\log\eps|
+\big|\!\log|\tau-\sigma|\big|
\big)
\,\rd\tau\,\rd\sigma\\
&\leq
c_d\,\eps^2\big(1+|\!\log\eps|\big).
\end{align*}
Finally, if $1<\alpha<d$, then
\[
\Ea(\chi_{T_\eps})
\leq
c_{\alpha,d}\,\eps^2
\int_{-1}^1
\int_{-1}^1
1\,\rd\tau\,\rd\sigma
\leq
c_{\alpha,d}\,\eps^2.
\]
This shows \eqref{eq:estim-energy}.
\end{proof}

\begin{lemma}\label{lem:thin-plate-upper-bound}
Let $d \in \N$, $d\geq2$ and $0<\alpha<d$. Then
\begin{equation}\label{eq:estim-determinant}
\Da(\chi_{T_\eps})
\leq
\begin{dcases}
c_{\alpha,d}\,\eps^{\alpha+1},
\qquad &0<\alpha<1,\\[1mm]
c_d\,\eps^2\big(1+|\!\log\eps|\big)^{\frac{d-1}{d}},
\qquad &\alpha=1,\\[1mm]
c_{\alpha,d}\,\eps^{\frac{2d+\alpha-1}{d}},
\qquad &1<\alpha<3,\\[1mm]
c_d\,\eps^{\frac{2d+2}{d}}
\big(1+|\!\log\eps|\big)^{\frac1d},
\qquad &\alpha=3,\\[1mm]
c_{\alpha,d}\,\eps^{\frac{2d+2}{d}},
\qquad &3<\alpha<d.
\end{dcases}
\end{equation}
\end{lemma}

\begin{proof}
We write $y=(\eta,\zeta)\in\R^{d-1}\times\R$.
If both $x+(\theta-1)y$ and
$x+\theta y$ belong to $T_\eps$,
then
\[
y
=
\big(x+\theta y\big)-\big(x+(\theta-1)y\big)
\in
T_\eps-T_\eps
\subseteq
B_2^{d-1}\times(-2\eps,2\eps).
\]
It follows that
\[
\Sa(\chi_{T_\eps})(x)
\preceq
\frac12
\int_{B_2^{d-1}}\int_{-2\eps}^{2\eps}
\big(|\eta|^2+\zeta^2\big)^{\frac{\alpha-d-2}{2}}
(\eta,\zeta)\otimes(\eta,\zeta)
\,\rd\zeta\,\rd\eta
=:
M_\eps.
\]

We now estimate $M_\eps$, which bounds $\Sa(\chi_{T_\eps})(x)$
uniformly in $x$. The integration domain $B_2^{d-1}\times(-2\eps,2\eps)$ and the weight
$\big(|\eta|^2+\zeta^2\big)^{(\alpha-d-2)/2}$ are invariant under reflection of any coordinate.
For $i\neq j$, reflecting the $i$-th coordinate changes the sign of the integrand defining $(M_\eps)_{ij}$, so all off-diagonal entries
vanish. By symmetry under permutations of the coordinates of $\eta$, the first $d-1$ diagonal entries are equal. Therefore, $M_\eps$ has the form
\[
M_\eps
=
\begin{pmatrix}
\lambda_\eps I_{d-1}&0\\
0&\mu_\eps
\end{pmatrix}\!,
\]
with
\[
\lambda_\eps
=
\frac{1}{2(d-1)}I_\eps^\parallel,
\qquad
\mu_\eps
=
\frac12I_\eps^\perp,
\]
where
\begin{align}
I_\eps^\parallel
&:=
\int_{B_2^{d-1}}\int_{-2\eps}^{2\eps}
\big(|\eta|^2+\zeta^2\big)^{\frac{\alpha-d-2}{2}}
|\eta|^2
\,\rd\zeta\,\rd\eta,
\label{Iparallel}
\\
I_\eps^\perp
&:=
\int_{B_2^{d-1}}\int_{-2\eps}^{2\eps}
\big(|\eta|^2+\zeta^2\big)^{\frac{\alpha-d-2}{2}}
\zeta^2
\,\rd\zeta\,\rd\eta.
\label{Iperp}
\end{align}
Thus,
\[
\det M_\eps=\lambda_\eps^{d-1}\mu_\eps.
\]

We proceed to estimate the quantities $I_\eps^\parallel$ and $I_\eps^\perp$.  
We start with $I_\eps^\parallel$ in \eqref{Iparallel} and split the domain of integration into two parts. On the region $|\eta|<2\eps$, we use
\[
\big(|\eta|^2+\zeta^2\big)^{\frac{\alpha-d-2}{2}}
|\eta|^2
\leq
|(\eta,\zeta)|^{\alpha-d},
\]
so the corresponding contribution is bounded by
$c_{\alpha,d} \, \eps^\alpha$. On the complementary region
$2\eps\leq|\eta|<2$, we use
\[
\big(|\eta|^2+\zeta^2\big)^{\frac{\alpha-d-2}{2}}
|\eta|^2
\leq
|\eta|^{\alpha-d}
\]
and its contribution is bounded by
\[
c_d \, \eps\int_{2\eps}^2r^{\alpha-2}\,\rd r.
\]
Combining the two contributions, we arrive at
\begin{equation}\label{eq:slab-parallel-estimate}
I_\eps^\parallel
\leq
\begin{dcases}
c_{\alpha,d} \, \eps^\alpha,
\qquad &0<\alpha<1,\\[1mm]
c_d \, \eps\big(1+|\!\log\eps|\big),
\qquad &\alpha=1,\\[1mm]
c_{\alpha,d} \, \eps,
\qquad &1<\alpha<d.
\end{dcases}
\end{equation}

Similarly, the contribution of $|\eta|<2\eps$ to $I_\eps^\perp$ is
bounded by $c_{\alpha,d} \, \eps^\alpha$. On the region
$2\eps\leq|\eta|<2$, we use
\[
\big(|\eta|^2+\zeta^2\big)^{\frac{\alpha-d-2}{2}}
\leq
|\eta|^{\alpha-d-2}
\]
and obtain the bound
\[
c_d \, \eps^3\int_{2\eps}^2r^{\alpha-4}\,\rd r.
\]
Again combining the two, we get
\begin{equation}\label{eq:slab-perpendicular-estimate}
I_\eps^\perp
\leq
\begin{dcases}
c_{\alpha,d} \, \eps^\alpha,
\qquad &0<\alpha<3,\\[1mm]
c_d \, \eps^3\big(1+|\!\log\eps|\big),
\qquad &\alpha=3,\\[1mm]
c_{\alpha,d} \, \eps^3,
\qquad &3<\alpha<d.
\end{dcases}
\end{equation}

By determinant monotonicity,
\[
\big(\!\det\Sa(\chi_{T_\eps})(x)\big)^{1/d}
\leq
\big(\!\det M_\eps\big)^{1/d}
=
c_d
\big(I_\eps^\parallel\big)^{\frac{d-1}{d}}
\big(I_\eps^\perp\big)^{\frac1d},
\]
and hence
\[
\Da(\chi_{T_\eps})
\leq
c_d \,\eps
\big(I_\eps^\parallel\big)^{\frac{d-1}{d}}
\big(I_\eps^\perp\big)^{\frac1d}.
\]
Then  \eqref{eq:estim-determinant} follows from
\eqref{eq:slab-parallel-estimate} and
\eqref{eq:slab-perpendicular-estimate}.
\end{proof}

\subsection{Proof of Theorem~\ref{thm:sharpness_1}}

\begin{proof}
We first prove part~{\rm (i)}. Fixing $r>p_\alpha$, set  $\rho_\eps:=\rho_{\eps,r}$ as defined in \eqref{eq:sharpness-normalized-thin-plate}. The density $\rho_\eps$ is log-concave, is supported in $B_{\sqrt2}$, and satisfies $\|\rho_\eps\|_{L^{r}(\R^d)}=1$ by \eqref{eq:thin-plate-lorentz} with $s=r$. By the homogeneity of $\Ea$ and
Lemma~\ref{lem:energy-thin-plate-upper-bound},
\[
\Ea(\rho_\eps)
\leq
\begin{dcases}
c_{\alpha,d,r} \, \eps^{\alpha+1-\frac2r},
\qquad &0<\alpha<1,\\[1mm]
c_{d,r} \, \eps^{2-\frac2r}\big(1+|\log\eps|\big),
\qquad &\alpha=1,\\[1mm]
c_{\alpha,d,r} \, \eps^{2-\frac2r},
\qquad &1<\alpha<d.
\end{dcases}
\]
If $0<\alpha<1$, then $r > p_\alpha$ is equivalent to
\[
\alpha+1-\frac2r>0.
\]
If $1\leq\alpha<d$, then $p_\alpha=1$, and $r>p_\alpha$ is equivalent
to $2-2/r>0$. The logarithmic factor at $\alpha=1$ is dominated by
the positive power of $\eps$. Hence
\[
\Ea(\rho_\eps)\to0
\qquad\text{as }\eps\to0.
\]
This proves part~{\rm (i)}.

We next prove part~{\rm (ii)}. We again set
$\rho_\eps:=\rho_{\eps,r}$, now with $r>q_\alpha$. Suppose first that $d=1$. Then $0<\alpha<1$, $\beta=3$, and
\[
q_\alpha=\frac3{\alpha+1}.
\]
The spatial scaling of $\Da$ gives
\[
\Da(\chi_{T_\eps})
=
\eps^{\alpha+1}\Da(\chi_{(-1,1)}).
\]
Since $|T_\eps|=2\eps$, the homogeneity of $\Da$ yields
\[
\Da(\rho_\eps)
=
c_{\alpha,r} \, \eps^{\alpha+1-\frac3r}.
\]
The exponent is positive precisely when $r>q_\alpha$, and therefore
$\Da(\rho_\eps)\to0$ as $\eps\to0$.

Assume next that $d\geq2$. The exponent of $\eps$ in Lemma~\ref{lem:thin-plate-upper-bound} is exactly $s_\alpha$ defined in~\eqref{eq:s_alpha}. Since
$|T_\eps|=2|B_1^{d-1}|\eps$, the homogeneity of $\Da$ gives
\begin{equation}\label{eq:sharpness-normalized-determinant-slab}
\Da(\rho_\eps)
\leq
c_{\alpha,d,r} \, \eps^{s_\alpha-\frac\beta r}
\begin{dcases}
\big(1+|\!\log\eps|\big)^{\frac{d-1}{d}},
\qquad &\alpha=1,\\[1mm]
\big(1+|\!\log\eps|\big)^{\frac1d},
\qquad &\alpha=3,\\[1mm]
1,
\qquad &\alpha\neq1,3.
\end{dcases}
\end{equation}
Since $q_\alpha=\beta/s_\alpha$, we have $s_\alpha-\frac\beta r>0$ 
whenever $r$ is selected to satisfy $r > q_\alpha$. The logarithmic factors in~\eqref{eq:sharpness-normalized-determinant-slab} are again dominated by the positive power of $\eps$. Hence
\[
\Da(\rho_\eps)\to0
\qquad\text{as }\eps\to0.
\]
This proves part~{\rm (ii)}.
\end{proof}

\begin{remark*}
The same examples establish the optimality of the first Lorentz indices in Theorem~\ref{thm:ED-quasiconcave}. Indeed, by~\eqref{eq:thin-plate-lorentz}, every Lorentz quasi-norm of $\rho_{\eps,r}$ with first index $r$ equals one, while, as $\eps\to0$, $\Ea(\rho_{\eps,r})\to0$ if $r>p_\alpha$ and $\Da(\rho_{\eps,r})\to0$ if $r>q_\alpha$.
Since these densities are log-concave and hence quasiconcave, neither first index can be increased, regardless of the choice of the second Lorentz index. A dilation argument then gives the same conclusion for densities supported in any prescribed ball $B_R$.
\end{remark*}

\section*{Acknowledgments}
The authors acknowledge the use of AI tools to assist with some computations in this work.
The authors developed the mathematical ideas and take full responsibility for the content and correctness of the manuscript.

The research of Alves was supported by King Abdullah University of Science and Technology (KAUST) under Award No.\ ORFS-CRG12-2024-6430. The research of Tzavaras was supported by KAUST baseline funds.

\end{document}